\documentclass[12pt,reqno]{amsart}

\usepackage{mathrsfs}
\usepackage{amsfonts}
\usepackage{amssymb}
\usepackage{amsthm}
\usepackage{amsmath}
\usepackage{latexsym}
\usepackage{tikz}
\usepackage{hyperref}
\usepackage{xcolor}

\hypersetup{hidelinks}
\usepackage{cite}
\allowdisplaybreaks

\newtheorem{theorem}{Theorem}[section]

\newtheorem{lemma}[theorem]{Lemma}

\newtheorem{remark}[theorem]{Remark}

\numberwithin{equation}{section}

\newcommand{\beq}{\begin{equation}}
\newcommand{\eeq}{\end{equation}}

\newcommand{\Rmnum}[1]{\expandafter\@slowromancap\romannumeral #1@}

\newcommand{\ben}{\begin{eqnarray}}
\newcommand{\een}{\end{eqnarray}}
\newcommand{\beno}{\begin{eqnarray*}}
\newcommand{\eeno}{\end{eqnarray*}}

\begin{document}

\title[{2D Anisotropic Navier--Stokes Equations}]
{Global Stability and Decay for the 2D Anisotropic Navier--Stokes Equations with Fractional Horizontal Dissipation}

\author[C. Y. Shen]{Chenyan Shen}
\address{Department of Mathematics, East China University of Science and Technology, Shanghai 200237,  P.R. China}
\email{18918598589@163.com}

\author[D. Wu]{Dan Wu}
\address{Department of Mathematics, East China University of Science and Technology, Shanghai 200237,  P.R. China}
\email{wudan202209@163.com}

\author[Y. Zhu]{Yi Zhu$^*$}
\address{Department of Mathematics, East China University of Science and Technology, Shanghai 200237,  P.R. China}
\email{zhuyim@ecust.edu.cn}
\thanks{$^*$ Corresponding author.}

\date{}

\begin{abstract}
	Although the classical Yudovich theory guarantees global existence and uniqueness of two-dimensional Euler solutions with bounded vorticity, their higher-order norms may still grow in time. This raises a fundamental question: how much dissipation we need to ensure global stability and decay? While the full Laplacian dissipation is sufficient, the stabilizing effect of weaker, anisotropic dissipation is considerably more subtle. In this paper, we investigate the global stability and long-time decay of the two-dimensional incompressible Navier–Stokes equations with fractional horizontal dissipation \(\Lambda_1^{2s}u\), where \(s<1\). A central challenge is to recover sufficient decay of the second velocity component, whose weak dissipation obstruct the closure of conventional anisotropic energy estimates. Previous work by Wang-Wu-Zhu [arXiv:2601.16110] resolved this difficulty whereas the spatial weighted energy functionals must be used when $s$ approaches $1$.
	The main purpose of our manuscript is to remove this restriction.
	 We introduce a hierarchy of carefully designed time-weighted anisotropic energy functionals which captures the interplay between horizontal fractional dissipation and the decay of velocity components. Under such initial assumptions, the characteristics of anisotropic dissipation is presented more straightforwardly, which also reveals the more essential intrinsic structure of fluid stability.
\end{abstract}

\maketitle

\section{introduction}

The Navier--Stokes equations with anisotropic dissipation arise naturally in fluid dynamics when horizontal and vertical length scales differ significantly.
In geophysical, thin-layer, and boundary-layer flows, suitable rescaling may lead to models in
which dissipation acts predominantly in one spatial direction. 
Among such anisotropic models, the two-dimensional incompressible Navier--Stokes equations with fractional dissipation in a single direction provide a prototypical example. In this paper, we consider the following Cauchy problem on $\mathbb{R}^2$ 
\begin{equation}\label{INS}
\begin{cases}
\partial_t u+\nu\Lambda_1^{2s}u+u\cdot\nabla u=-\nabla p, \qquad (t,x)\in \mathbb{R}^{+} \times\mathbb{R}^2
\\
\nabla\cdot u=0,
\\
u|_{t=0}=u_0,
\end{cases}
\end{equation}
where $u=(u_1,u_2)$ denotes the velocity field, $p$ is the pressure,
$\nu>0$ is the viscosity coefficient. The fractional operator $\Lambda_1^{2s}$
acts only in the horizontal direction and is
defined by
$$
\widehat{\Lambda_1^{2s}f}(\xi)
\triangleq
|\xi_1|^{2s}\widehat f(\xi),
\qquad 0\leq s\leq 1 .
$$
When $s=1$, system \eqref{INS} reduces to the Navier--Stokes equations
with horizontal Laplacian dissipation, while $0<s<1$ corresponds to a
fractional anisotropic dissipative model. Besides, when $s=0$, the dissipation term will reduce to a standard damping term.

Mathematically, system \eqref{INS} lies between the inviscid Euler equations and the classical fully dissipative Navier--Stokes equations. Recall that the incompressible Euler equations
\begin{equation*}
\begin{cases}
\partial_t u+u\cdot\nabla u=-\nabla p,
\\
\nabla\cdot u=0,
\end{cases}
\end{equation*}
describe the motion of ideal incompressible fluids without viscosity. 
Although classical Yudovich theory guarantees global existence and uniqueness of two-dimensional Euler solutions with bounded vorticity, higher-order Sobolev norms of smooth solutions may still grow
in time. For example, Kiselev and \v{S}ver\'ak \cite{KS} constructed
smooth solutions in a disk for which the gradient of vorticity grows at a
double exponential rate. 
Similar growth phenomena have also been established in other settings, including exponential growth on the torus \cite{Z}, linear growth in the whole space \cite{C}, and superlinear growth on $\mathbb{T}^2$ \cite{D1} and in $\mathbb{R}^2$ \cite{YZ};    see also \cite{D2,D3,D} for related results. 
These results show that the absence of dissipation may hinder the stabilization of solutions to the Euler equations. A natural question is how much dissipation is needed
to recover stability and decay? Full Laplacian dissipation is
certainly sufficient. As demonstrated by Schonbek and others (see, for example, \cite{S1,S2,S3,W}), solutions of the Navier–Stokes equations
\begin{equation*}
\begin{cases}
\partial_t u-\nu\Delta u+u\cdot\nabla u=-\nabla p,
\\
\nabla\cdot u=0,
\end{cases}
\end{equation*}
are stable and decay at explicit decay rates.

However, the analysis becomes much more delicate when the dissipation acts in 
only one spatial direction. In this anisotropic setting, there is no vertical dissipation, and the nonlinear term $u_2\partial_2 u$ is not directly controlled by the available
horizontal dissipation. This difficulty already appears in the classical
case $s=1$, and the global stability problem in the whole space $\mathbb R^2$
remains highly nontrivial.
Nevertheless, several
important results have been obtained in related settings. In the partially 
periodic setting $\mathbb{T}\times\mathbb{R}$, Dong, Wu, Xu and Zhu 
\cite{DWXZ} established global stability for the two-dimensional Navier–Stokes equations with horizontal Laplacian dissipation ($s=1$) and proved exponential decay of the oscillatory part for small initial data. Their proof relies on a decomposition of the velocity into its horizontal average and oscillatory part. 
In three dimensions, Ji, Wu 
and Yang \cite{J} obtained stability and optimal decay estimates for the 
Navier--Stokes equations with only horizontal dissipation. Recently, Bie, Fang, Wang and Zhou \cite{Bie} extended this result to the case of fractional horizontal dissipation, 
and established global stability and sharp decay estimates under a broader range of dissipation exponents. 
For further results on stability and large-time behavior for related
partially dissipative fluid models, especially anisotropic
magnetohydrodynamic systems, we refer to \cite{Zhu,Zhu1,Lin1,Lin2,zz,zz1,lai1,lai2,lai3,zhang,Du,Jiang,li,Ren}.

For the two-dimensional fractional problem, Wang, Wu, and Zhu
\cite{WWZ} established global stability and decay estimates for the
full range $0\leq s<1$. When $0\leq s\leq \frac34$, the high-order energy
estimates can be closed directly through a suitable redistribution of
horizontal derivatives. For $\frac34<s<\frac{11}{12}$, their proof combines
anisotropic energy estimates with large-time decay obtained from the
Duhamel formulation. 
In particular, the stream-function representation and incompressibility yield improved decay of the second velocity component $u_2$. To cover the remaining range $\frac{11}{12}\leq s<1$, they introduce polynomial weights in the vertical variable and use weighted $L^2$ estimates for Riesz transforms to obtain sufficient decay.

The purpose of the present paper is to provide an alternative treatment
of the regime $\frac34<s<1$ without physical-space weights. 
We explain below the main
difficulties and the role of decay estimates in the proof.

\begin{enumerate}
	
	\item[$\bullet$] \textit{The higher-order energy estimate.}
		We first seek a uniform bound for
		$$
		\sup_{0\leq\tau\leq t}\|u(\tau)\|_{H^k}^2
		+
		\int_0^t\|\Lambda_1^s u(\tau)\|_{H^k}^2\,\mathrm{d}\tau.
		$$
		The difficulty comes from the nonlinear terms involving
		the highest-order vertical derivatives.
		For example, the estimate for $\|\partial_2^k u\|_{L^2}$
		will lead to 
		\begin{equation}\label{iden}
			\mathcal N_k
			=
			\int_{\mathbb R^2}
			\partial_2u_2\,|\partial_2^k u_1|^2\,\mathrm{d}x
			=
			-\int_{\mathbb R^2}
			\partial_1u_1\,|\partial_2^k u_1|^2\,\mathrm{d}x,
		\end{equation}
		where we have used the divergence-free condition.
		This term \eqref{iden} contains one horizontal derivative, however
		the highest-order derivatives remain purely vertical. Anisotropic type estimates and interpolation then give 
		\begin{equation}\label{N_k}
			\begin{aligned}
				|\mathcal N_k|
				&\lesssim
				\|\Lambda_1^s u\|_{L^2_{x_1}L^\infty_{x_2}}
				\|
				\Lambda_1^{\frac34-\frac{s}{2}}\partial_2^k u
				\|_{L^2}^2
				\\
				&\lesssim
				\|u\|_{H^k}^{3-\frac{3}{2s}}
				\|\Lambda_1^s u\|_{H^k}^{\frac{3}{2s}}.
			\end{aligned}
		\end{equation}
	Such estimate \eqref{N_k} does not close the global uniform 
		$H^k$ bound since $\frac{3}{2s}<2$. To overcome this difficulty, we have to refine the
		nonlinear estimates by exploiting the decay of 
		individual velocity components and their derivatives.

	\item[$\bullet$] \textit{Time decay of lower-order quantities.}
		We impose the additional low-frequency
		condition
		$$
		\|\Lambda_1^{-\sigma_1}u_0\|_{L^2}\ll1,
		\qquad
		\frac{s}{2}<\sigma_1<\frac12.
		$$
		The negative-order norm's estimate for $u$ will lead to 
		$u_2\partial_2u$, while the time-weighted estimate for
		$\partial_1u$ involves $\partial_1u_2\partial_2u$.
		To control these terms, we have to discover additional decay property for $u_2$.
		The vorticity formulation and the Biot--Savart law yield {(with $\nabla^\perp=(-\partial_2,\partial_1)$ and
			$\nabla^\perp \cdot F=\partial_1F_2-\partial_2F_1$)}
		$$
		\partial_tu_2+\nu\Lambda_1^{2s}u_2
		=
		(-\Delta)^{-1}\partial_1
		\nabla^\perp\!\cdot(u\cdot\nabla u).
		$$
		The explicit horizontal derivative allows one additional
		inverse horizontal derivative and
		we therefore also assume that 
		$$
		\|\Lambda_1^{-\sigma_1-1}\Lambda^\sigma u_{2,0}\|_{L^2}
		\ll1,
		\qquad
		\sigma_1<\sigma<\frac12
		.
		$$
		The factor $\Lambda^\sigma$ weakens the low-frequency
		singularity in the nonlinear estimates.
		Using these negative-order bounds in the time-weighted
		energy frame we shall derive 
		$$
		\begin{aligned}
			\|u(t)\|_{L^2}
			&\lesssim (1+t)^{-\frac{\sigma_1}{2s}},\\
			\|\partial_1u(t)\|_{L^2}
			&\lesssim (1+t)^{-\frac{\sigma_1+1}{2s}},\\
			\|\Lambda^\sigma u_2(t)\|_{L^2}
			&\lesssim (1+t)^{-\frac{\sigma_1+1}{2s}}.
		\end{aligned}
		$$
		The resulting decay properties supply the time
		integrability needed in the higher-order estimate.
\end{enumerate}

Our main result can be stated as follows.

\begin{theorem}\label{thm1.1}
Consider the anisotropic Navier--Stokes system \eqref{INS} with viscosity coefficient $\nu>0$. Let $s,\sigma_1,\sigma$ satisfy
$$
\frac34<s<2\sigma_1,\qquad\sigma_1<\sigma<\frac12,
$$
and let $k\in\mathbb N$ satisfy
$$
k>\max\left\{
\frac{\frac{11}{4}+7\sigma_1-4\sigma+\sigma_1(2\sigma_1+\sigma)}{1-2\sigma_1},\frac{3\sigma-\sigma_1-4\sigma_1^2}
{2(\sigma-\sigma_1)}\right\},\qquad k\in\mathbb N.
$$
Then there exists a sufficiently small constant $\epsilon>0$ such that, for any divergence-free initial data $u_0\in H^k(\mathbb R^2)$ satisfying 
\begin{equation}\label{initial}
\|u_0\|_{H^k}+\|\Lambda^{-\sigma_1}_1u_0\|_{L^2}+
\|\Lambda^{-\sigma_1-1}_1\Lambda^{\sigma} u_{2,0}\|_{L^2}\leq  \epsilon,
\end{equation}
system \eqref{INS} admits a unique global-in-time solution
$$
u\in C([0,\infty);H^k(\mathbb{R}^2)),
$$
satisfying
$$
\|u\|^2_{H^k}+\int_0^t\|\Lambda^s_1u(\tau)\|^2_{H^k}\,\mathrm{d}\tau\lesssim \epsilon^2.
$$
Furthermore, the solution $u$ admits the following upper bounds, for any $t>0$, 
$$
\sup_{0\le \tau\le t}
\|\Lambda_1^{-\sigma_1}u(\tau)\|_{L^{2}}^2
+
\int_0^t
\|\Lambda_1^{s-\sigma_1}u(\tau)\|_{L^{2}}^2\,\mathrm{d}\tau
\lesssim \epsilon^2,
$$
$$
\sup_{0\le \tau\le t}
(1+\tau)^{\frac{\sigma_1}{s}}
\|u(\tau)\|_{L^{2}}^2
+
\int_0^t
(1+\tau)^{\frac{\sigma_1}{s}}
\|\Lambda_1^s u(\tau)\|_{L^{2}}^2\,\mathrm{d}\tau \lesssim \epsilon^2,
$$
$$
\sup_{0\le \tau\le t}
(1+\tau)^{\frac{\sigma_1+1}{s}}
\|\partial_1u(\tau)\|_{ L^{2}}^2
+
\int_0^t
(1+\tau)^{\frac{\sigma_1+1}{s}}
\|\Lambda_1^{1+s}u(\tau)\|_{ L^{2}}^2\,\mathrm{d}\tau  \lesssim \epsilon^2,
$$
$$
\sup_{0\le \tau\le t}
\|\Lambda_1^{-\sigma_1-1}\Lambda^\sigma u_2(\tau)\|_{L^2}^2
+
\int_0^t
\|\Lambda_1^{-\sigma_1-1+s}\Lambda^\sigma u_2(\tau)\|_{L^2}^2\,\mathrm{d}\tau  \lesssim \epsilon^2,
$$
$$
\sup_{0\le \tau\le t}
(1+\tau)^{\frac{\sigma_1+1}{s}}
\|\Lambda^\sigma u_2(\tau)\|_{L^2}^2
+
\int_0^t
(1+\tau)^{\frac{\sigma_1+1}{s}}
\|\Lambda_1^s\Lambda^\sigma u_2(\tau)\|_{L^2}^2\,\mathrm{d}\tau \lesssim \epsilon^2.
$$

\end{theorem}

\begin{remark}
For the case $0\leq s\leq \frac34$, the stability and decay theorem can
be obtained directly from traditional energy method, and we
omit the details here. Therefore, we just focus on the more difficult regime
$\frac34<s<1$.
\end{remark}

\paragraph{\textbf{Notations.}}

{We use the Fourier transform on $\mathbb R^2$,}
\[
\widehat f(\xi)=\frac1{2\pi}\int_{\mathbb R^2}e^{-ix\cdot\xi}f(x)\,\mathrm{d}x,
\qquad
\check h(x)=\mathcal F^{-1}h(x)
=\frac1{2\pi}\int_{\mathbb R^2}e^{ix\cdot\xi}h(\xi)\,\mathrm{d}\xi.
\]
For $l\in\mathbb R$, define $\Lambda^l$ and
		$\Lambda_j^l$ $(j=1,2)$ by
		$$
		\widehat{\Lambda^l f}(\xi)
		=|\xi|^l\widehat f(\xi),
		\qquad
		\widehat{\Lambda_j^l f}(\xi)
		=|\xi_j|^l\widehat f(\xi).
		$$
		In particular, $\Lambda=(-\Delta)^{\frac12}$ and
		$\Lambda_j=(-\partial_j^2)^{\frac12}$.
For $f\in L^2$, define the auxiliary function as follows
\[
f^\sharp \triangleq \check{\lvert\widehat f\rvert}
=\mathcal F^{-1}(\lvert\widehat f\rvert).
\]
This function has the Fourier transform $\lvert\widehat f\rvert$, 
it is generally different from $f$ and not need to be nonnegative in 
physical space.

 Throughout this paper, $C>0$ stands for a generic constant which may change line by line. For brevity, we write
$f\lesssim g$ instead of $f\leq Cg$. The notation $f\approx g$ means that $
f\lesssim g$ and $g\lesssim f$. Also, $L^p(\mathbb{R}^d)$ denotes the standard
Lebesgue space on the $d$-dimensional Euclidean space $\mathbb{R}^d$.
For $f,g\in L^2(\mathbb{R}^d)$, we denote the $L^2$ inner product by
$$
\langle f,g\rangle \triangleq\int_{\mathbb{R}^d}f(x)g(x)\,\mathrm{d}x.
$$
We also use anisotropic spaces. Let $N_1$ and $N_2$ be two normed spaces
defined on the measure spaces $(X_1,\mu_1)$ and $(X_2,\mu_2)$,
respectively. For a function $f$ defined on the product space
$X_1\times X_2$, the anisotropic space
$N_1(X_1,\mu_1)N_2(X_2,\mu_2)$ is defined by
$$
\|f\|_{N_1(X_1,\mu_1)N_2(X_2,\mu_2)}
:=
\left\|
\|f(x_1,x_2)\|_{N_1(X_1,\mu_1)}
\right\|_{N_2(X_2, \mu_2)}.
$$
In particular, we use the notation $L^p_{x_1}L^q_{x_2}$ to denote the
anisotropic Lebesgue space obtained by taking the $L^p$-norm with respect
to the $x_1$ variable first and then the $L^q$-norm with respect to the
$x_2$ variable. For $k\geq0$, the Sobolev space $H^k(\mathbb{R}^d)$ is defined by
$$
\|f\|_{H^k}
\triangleq
\|(1+|\xi|^2)^{\frac k2}\widehat f(\xi)\|_{L^2_\xi},
$$
while the homogeneous Sobolev space $\dot H^k(\mathbb{R}^d)$ is equipped
with the norm
$$
\|f\|_{\dot H^k}
\triangleq
\||\xi|^k\widehat f(\xi)\|_{L^2_\xi}.
$$

The paper is organized as follows. 
We collect the preliminary lemmas in Section 2
and prove Theorem 1.1 in Section 3.

\section{Preliminaries}\label{pre}

The fractional dissipation in the system acts only in the horizontal direction, which leads to an anisotropic structure in the nonlinear estimates. We therefore begin with several anisotropic product inequalities that will be used throughout the paper. These estimates are standard tools in the analysis of anisotropic partial differential equations, see, for example, \cite{Cao}.

\begin{lemma}\label{2.1}
Let $f,g\in H^1(\mathbb{R}^2).$ Then the following estimates hold:
\begin{equation*}
\|fg\|_{L^1_{x_1}L^2_{x_2}}\lesssim \|f\|^{\frac12}_{L^2}\|\partial_2f\|^{\frac12}_{L^2}\|g\|_{L^2},
\end{equation*}
and
\begin{equation*}
\|fg\|_{L^1_{x_2}L^2_{x_1}}\lesssim \|f\|^{\frac12}_{L^2}\|\partial_1f\|^{\frac12}_{L^2}\|g\|_{L^2}.
\end{equation*}
\end{lemma}
\begin{proof}
By Hölder's inequality in the $x_1$ variable, we have
$$
\|fg\|_{L^1_{x_1}L^2_{x_2}}
=
\left\|
\|fg\|_{L^1_{x_1}}
\right\|_{L^2_{x_2}}
\lesssim
\left\|
\|f\|_{L^2_{x_1}}
\|g\|_{L^2_{x_1}}
\right\|_{L^2_{x_2}}
\lesssim
\left\|
\|f\|_{L^2_{x_1}}
\right\|_{L^{\infty}_{x_2}}
\|g\|_{L^2} .
$$
Applying Minkowski's inequality and the one-dimensional
Gagliardo--Nirenberg inequality in the $x_2$ variable yields
$$
\left\|
\|f\|_{L^2_{x_1}}
\right\|_{L^\infty_{x_2}}
\lesssim
\|f\|_{L^2}^{\frac12}
\|\partial_2 f\|_{L^2}^{\frac12}.
$$
Combining the above estimates yields the first inequality and the second one  follows similarly.
\end{proof}
The following classical Hardy--Littlewood--Sobolev inequality will be
used repeatedly in the subsequent estimates and the proof can be
found in \cite{Dan}.
\begin{lemma}[Hardy--Littlewood--Sobolev inequality]\label{2.2}
The following inequality holds
\begin{equation*}
\|\Lambda^{-\rho}f\|_{L^q(\mathbb{R}^d)}
\lesssim\|f\|_{L^p(\mathbb{R}^d)},
\end{equation*}
with
$$
\frac{1}{q}+\frac{\rho}{d}=\frac{1}{p},
\qquad
0<\rho<d,
\qquad
1<p<q<\infty .
$$
\end{lemma}
We shall use the following two classical product estimates. 
The first is the standard product law in homogeneous Sobolev spaces, 
see, for instance, \cite{Ada}.

\begin{lemma}\label{2.3}
	Let $f$ and $g$ be sufficiently smooth functions on $\mathbb{R}^d$.
	Assume that
	$$
	s_1,s_2<\frac{d}{2},\qquad s_1+s_2>0,
	$$
	and
	$$
	s+\frac{d}{2}=s_1+s_2.
	$$
	Then there holds 
	$$
	\|\Lambda^s(fg)\|_{L^2}
	\lesssim
	\|\Lambda^{s_1}f\|_{L^2}
	\|\Lambda^{s_2}g\|_{L^2}.
	$$
\end{lemma}
The second is the classical fractional Leibniz inequality, see, for
example, \cite{LLL}.

\begin{lemma}\label{lem2.4}
	Let $m>0$ and $1<r<\infty$. Assume that
	$1<p_1,p_2,q_1,q_2\leq\infty$ satisfy
	$$
	\frac1r
	=
	\frac1{p_1}+\frac1{q_1}
	=
	\frac1{p_2}+\frac1{q_2}.
	$$
	Then, for sufficiently smooth functions $f$ and $g$ on $\mathbb{R}^d$, there holds 
	$$
	\|\Lambda^m(fg)\|_{L^r}
	\lesssim
	\|\Lambda^m f\|_{L^{p_1}}\|g\|_{L^{q_1}}
	+
	\|f\|_{L^{p_2}}\|\Lambda^m g\|_{L^{q_2}}.
	$$
\end{lemma}

Based on these classical product estimates, we shall  derive the following
anisotropic fractional product {estimate} tailored to the nonlinear
terms we need.

\begin{lemma}\label{lem:negative-product-1}
Let $0<a<\frac{1}{2}$ and $b>0$. For sufficiently smooth functions $f$ and $g$ on $\mathbb{R}^2$,
\[
\begin{split}
\|\Lambda_1^{-a}\Lambda^b(fg)\|_{L^2}
\lesssim{}&\|\Lambda^b f^\sharp\|_{L^2}
\|\Lambda_1^{\frac12-a}g^\sharp\|_{L^2_{x_1}L^\infty_{x_2}}\\
&+\|\Lambda_1^{\frac12-a}f^\sharp\|_{L^2_{x_1}L^\infty_{x_2}}
\|\Lambda^b g^\sharp\|_{L^2}.
\end{split}
\]
The implicit constant depends only on $a$ and $b$.
\end{lemma}
\begin{proof}
The convolution identity for the Fourier transform and the inequality
$|\xi|^b\lesssim|\xi-\eta|^b+|\eta|^b$ give
\[
\begin{split}
|\xi|^b|\widehat{fg}(\xi)|
\lesssim{}&\int_{\mathbb R^2}|\xi-\eta|^b
|\widehat f(\xi-\eta)|\,|\widehat g(\eta)|\, \mathrm{d}\eta\\
&+\int_{\mathbb R^2}|\widehat f(\xi-\eta)|
|\eta|^b|\widehat g(\eta)|\,\mathrm{d}\eta.
\end{split}
\]
Since these Fourier transforms are nonnegative, Plancherel's theorem implies, 
\begin{equation}\label{eq:fourier-envelope-splitting}
\begin{split}
\|\Lambda_1^r\Lambda^b(fg)\|_{L^2}
\lesssim{}&\|\Lambda_1^r((\Lambda^b f^\sharp)g^\sharp)\|_{L^2}\\
&+\|\Lambda_1^r(f^\sharp(\Lambda^b g^\sharp))\|_{L^2},
\end{split}
\end{equation}
for any $r\in\mathbb R$ when the right-hand side is bounded. 

Set $p_a=\frac{2}{1+2a}$ and $q_a=\frac{1}{a}$, then
$1<p_a<2<q_a<\infty$ and $\frac{1}{p_a}=\frac12+\frac{1}{q_a}.$
For each fixed $x_2$, we apply Lemma~\ref{2.2}
	in $x_1$ and then use H\"older's inequality
	in both variables, there is 
\[
\begin{split}
\|\Lambda_1^{-a}((\Lambda^b f^\sharp)g^\sharp)\|_{L^2}
&\lesssim\|(\Lambda^b f^\sharp)g^\sharp\|_{L^{p_a}_{x_1}L^2_{x_2}}\\
&\le\|\Lambda^b f^\sharp\|_{L^2}
\|g^\sharp\|_{L^{q_a}_{x_1}L^\infty_{x_2}}\\
&\lesssim\|\Lambda^b f^\sharp\|_{L^2}
\|\Lambda_1^{\frac12-a}g^\sharp\|_{L^2_{x_1}L^\infty_{x_2}}.
\end{split}
\]
The last inequality follows from the one-dimensional homogeneous Sobolev embedding
$\dot H^{\frac12-a}(\mathbb R)\hookrightarrow L^{\frac1a}(\mathbb R)$. The other product is symmetric.
We use \eqref{eq:fourier-envelope-splitting} with $r=-a$ to finish the proof. 
\end{proof}

\begin{lemma}\label{lem:fractoinal-product-1}
Let $\frac{3}{4}<s<1$, $\sigma>0$, and $d_s=\frac{3}{4}-\frac{s}{2}$.
For sufficiently smooth functions $f$ and $g$ on $\mathbb{R}^2$, there holds 
\[
\begin{split}
\|\Lambda_1^{1-s}\Lambda^\sigma(fg)\|_{L^2}
\lesssim{}&\|\Lambda_1^{d_s}f^\sharp\|_{L^2_{x_1}L^\infty_{x_2}}
\|\Lambda_1^{d_s}\Lambda^\sigma g^\sharp\|_{L^2}\\
&+\|\Lambda_1^{d_s}g^\sharp\|_{L^2_{x_1}L^\infty_{x_2}}
\|\Lambda_1^{d_s}\Lambda^\sigma f^\sharp\|_{L^2}.
\end{split}
\]
The implicit constant depends only on $s$ and $\sigma$.
\end{lemma}
\begin{proof}
Use \eqref{eq:fourier-envelope-splitting} with $b=\sigma$ and $r=1-s$, we still need to 
estimate
$\Lambda_1^{1-s}((\Lambda^\sigma f^\sharp)g^\sharp)$ and its
symmetric counterpart. Since
\[
\frac14<d_s<\frac38<\frac12,\qquad
2d_s>0,\qquad 2d_s-\frac12=1-s,
\]
we apply Lemma~\ref{2.3} in one dimensional case with
	$s_1=s_2=d_s$. 
	For each fixed $x_2$ it yields that, 
\[
\begin{split}
\|\Lambda_1^{1-s}((\Lambda^\sigma f^\sharp)g^\sharp)
(\cdot,x_2)\|_{L^2_{x_1}}
\lesssim
\|\Lambda_1^{d_s}\Lambda^\sigma f^\sharp(\cdot,x_2)\|_{L^2_{x_1}}
\|\Lambda_1^{d_s}g^\sharp(\cdot,x_2)\|_{L^2_{x_1}}.
\end{split}
\]
Take the $L^2$ norm in $x_2$ and put the second factor in
$L^\infty_{x_2}$ we shall obtain 
\[
\begin{split}
&\|\Lambda_1^{1-s}((\Lambda^\sigma f^\sharp)g^\sharp)\|_{L^2}\\
&\qquad\lesssim
\|\Lambda_1^{d_s}\Lambda^\sigma f^\sharp\|_{L^2}
\|\Lambda_1^{d_s}g^\sharp\|_{L^2_{x_1}L^\infty_{x_2}}.
\end{split}
\]
We then interchange $f$ and $g$  to estimate the other terms. Adding the two bounds will  complete the proof.
\end{proof}

\section{Energy Estimate and the Proof of Main Result}\label{proof}
\subsection{Construction of the energy framework}

We now introduce a hierarchy of time-weighted energy functionals, which will
be used to establish the global-in-time stability and decay estimates.
Motivated by the anisotropic dissipation mechanism of the system \eqref{INS},
we define the following energy functionals.

Throughout this section, we assume that
$$
\frac34<s<2\sigma_1,\qquad \sigma_1<\sigma<\frac12.
$$
We choose the Sobolev regularity index $k$ sufficiently large depending on $\sigma_1$ and $\sigma$. The precise condition on $k$ will be specified later. Now, we set 
\begin{equation}\label{eq1}
\begin{aligned}
E_0(t)
&=
\sup_{0\le \tau\le t}
\|\Lambda_1^{-\sigma_1}u(\tau)\|_{L^{2}}^2
+
\int_0^t
\|\Lambda_1^{s-\sigma_1}u(\tau)\|_{L^{2}}^2\,\mathrm{d}\tau,
\\[0.4em]
E_1(t)
&=
\sup_{0\le \tau\le t}
(1+\tau)^{\frac{\sigma_1}{s}}
\|u(\tau)\|_{L^{2}}^2
+
\int_0^t
(1+\tau)^{\frac{\sigma_1}{s}}
\|\Lambda_1^s u(\tau)\|_{L^{2}}^2\,\mathrm{d}\tau,
\\[0.4em]
E_2(t)
&=
\sup_{0\le \tau\le t}
\|u(\tau)\|_{ H^k}^2
+
\int_0^t
\|\Lambda_1^s u(\tau)\|_{ H^k}^2\,\mathrm{d}\tau,
\\[0.4em]
E_3(t)
&=
\sup_{0\le \tau\le t}
(1+\tau)^{\frac{\sigma_1+1}{s}}
\|\partial_1u(\tau)\|_{ L^{2}}^2
+
\int_0^t
(1+\tau)^{\frac{\sigma_1+1}{s}}
\|\Lambda_1^{1+s}u(\tau)\|_{ L^{2}}^2\,\mathrm{d}\tau,
\\[0.4em]
E_4(t)
&=
\sup_{0\le \tau\le t}
\|\Lambda_1^{-\sigma_1-1}\Lambda^\sigma u_2(\tau)\|_{L^2}^2
+
\int_0^t
\|\Lambda_1^{-\sigma_1-1+s}\Lambda^\sigma u_2(\tau)\|_{L^2}^2\,\mathrm{d}\tau,
\\[0.4em]
E_5(t)
&=
\sup_{0\le \tau\le t}
(1+\tau)^{\frac{\sigma_1+1}{s}}
\|\Lambda^\sigma u_2(\tau)\|_{L^2}^2
+
\int_0^t
(1+\tau)^{\frac{\sigma_1+1}{s}}
\|\Lambda_1^s\Lambda^\sigma u_2(\tau)\|_{L^2}^2\,\mathrm{d}\tau.
\\[0.4em]
\end{aligned}
\end{equation}
The energy $E_1(t)$ is an auxiliary quantity and can be controlled by
$E_0(t)$ and $E_3(t)$. Indeed, interpolation in the $x_1$ variable gives
$$
\|u(\tau)\|_{L^2}
\lesssim
\|\Lambda_1^{-\sigma_1}u(\tau)\|_{L^2}^{\frac1{1+\sigma_1}}
\|\partial_1u(\tau)\|_{L^2}^{\frac{\sigma_1}{1+\sigma_1}},
$$
and
$$
\|\Lambda_1^s u(\tau)\|_{L^2}
\lesssim
\|\Lambda_1^{s-\sigma_1}u(\tau)\|_{L^2}^{\frac1{1+\sigma_1}}
\|\Lambda_1^{1+s}u(\tau)\|_{L^2}^{\frac{\sigma_1}{1+\sigma_1}}.
$$
Consequently, 
$$
	\begin{aligned}
		&\sup_{0\leq\tau\leq t}
		(1+\tau)^{\frac{\sigma_1}{s}}\|u(\tau)\|_{L^2}^2\lesssim
		E^{\frac1{1+\sigma_1}}_0(t)
		E^{\frac{\sigma_1}{1+\sigma_1}}_3(t).
	\end{aligned}
$$
Moreover, by H\"older's inequality in time there is 
$$
\begin{aligned}
&\int_0^t
(1+\tau)^{\frac{\sigma_1}{s}}
\|\Lambda_1^s u(\tau)\|_{L^2}^2\,\mathrm{d}\tau\\&
\lesssim
\left(
\int_0^t
\|\Lambda_1^{s-\sigma_1}u(\tau)\|_{L^2}^2
\,\mathrm{d}\tau
\right)^{\frac1{1+\sigma_1}}
\left(
\int_0^t
(1+\tau)^{\frac{1+\sigma_1}{s}}
\|\Lambda_1^{1+s}u(\tau)\|_{L^2}^2
\,\mathrm{d}\tau
\right)^{\frac{\sigma_1}{1+\sigma_1}}.
\end{aligned}
$$
	Combining the inequalities above and applying Young's inequality, 
	we obtain
	\begin{equation}\label{E1-interpolation}
		E_1(t)
		\lesssim
		E_0(t)+E_3(t).
	\end{equation}
Thus, $E_1(t)$ can be excluded from the final energy frame.

\subsection{$\textit{A Priori}$ Estimates}
We now turn to the estimate of energy functionals introduced above. 
We first focus on the negative-order energy $E_0(t)$.
\begin{lemma}\label{lem1}
Assume that $$k>\max\left\{
\frac{6\sigma_1^2-6\sigma_1+2}{3-6\sigma_1},
\frac{\left|-6\sigma_1^2+6\sigma\sigma_1+4\sigma_1-3\right|}
{3-6\sigma_1},
1
\right\},$$
and the energies are defined as in \eqref{eq1}, then we have
$$
E_0(t)\lesssim  E_0(0)+ E^{\frac32}_0(t)+E^{\frac32}_1(t)+E^{\frac32}_2(t)+E^{\frac32}_3(t)+E^{\frac32}_5(t).
$$
\end{lemma}
\begin{proof}
Applying the operator $\Lambda_1^{-\sigma_1}$ to system \eqref{INS} and taking  $L^2$ inner product with $\Lambda_1^{-\sigma_1}u$, we obtain
\begin{equation*}
\begin{split}
\frac12\frac{\mathrm{d}}{\mathrm{d}t}\|\Lambda_1^{-\sigma_1}u\|^2_{L^2}+\nu\|\Lambda_1^{s-\sigma_1}u\|^2_{L^2}&=-<\Lambda_1^{-\sigma_1}(u\cdot\nabla u),\Lambda_1^{-\sigma_1}u>\\&=-<\Lambda_1^{-\sigma_1}(u_1\partial_1u),\Lambda_1^{-\sigma_1}u>\\&\quad\,\,-<\Lambda_1^{-\sigma_1}(u_2\partial_2 u),\Lambda_1^{-\sigma_1}u>\\&\triangleq J_1+J_2.
\end{split}
\end{equation*}
{\bf{Estimate of $J_1$:}} 
Applying Lemma \ref{2.2} in the $x_1$ variable and using H\"older's inequality,
we shall obtain
\begin{equation*}
\begin{split}
|J_1|&\lesssim
\|\Lambda_1^{-\sigma_1}(u_1\partial_1u)\|_{L^2}\|\Lambda_1^{-\sigma_1}u\|_{L^2}\\&\lesssim \|u_1\partial_1u\|_{L^{\frac{2}{1+2\sigma_1}}_{x_1}L^2_{x_2}}\|\Lambda_1^{-\sigma_1}u\|_{L^2}\\&\lesssim 
\|u_1\|_{L^{\frac{1}{\sigma_1}}_{x_1}L^{\infty}_{x_2}}\|\partial_1u\|_{L^2}\|\Lambda_1^{-\sigma_1}u\|_{L^2}.
\end{split}
\end{equation*}
Anisotropic type Sobolev embedding theorem and interpolation yield
\begin{equation}\label{J_11}
\begin{aligned}
\|u_1\|_{L_{x_1}^{\frac1{\sigma_1}}L_{x_2}^{\infty}}
&\lesssim
\|\Lambda_1^{\frac12-\sigma_1}u_1\|_{L_{x_1}^2L_{x_2}^{\infty}}
\\
&\lesssim
\|\Lambda_1^{\frac12-\sigma_1}u_1\|_{L^2}^{\theta_1}
\|\Lambda_1^{\frac12-\sigma_1}
\Lambda_2^{k-\frac12+\sigma_1}u_1\|_{L^2}^{1-\theta_1}
\\
&\lesssim
\|u_1\|_{L^2}^{\theta_1\theta_2}
\|\Lambda_1^s u_1\|_{L^2}^{\theta_1(1-\theta_2)}
\|u_1\|_{H^k}^{1-\theta_1}.
\end{aligned}
\end{equation}
Then H\"older's inequality will give 
\begin{equation}\label{J_1}
\begin{split}
\int_0^{t}|J_1|\,\mathrm{d}\tau& \lesssim \sup_{0\le \tau\le t}
(1+\tau)^{\frac{\sigma_1\theta_1\theta_2}{2s}}\|u_1\|^{\theta_1\theta_2}_{L^2}\left(\int_0^t(1+\tau)^{\frac{\sigma_1}{s}}\|\Lambda^s_1u_1\|^2\,\mathrm{d}\tau\right)^{\frac{\theta_1(1-\theta_2)}{2}}\\&\quad\,\,\cdot
\sup_{0\le \tau\le t}\|u_1\|^{1-\theta_1}_{H^k}\sup_{0\le \tau\le t}
(1+\tau)^{\frac{\sigma_1+1}{2s}}\|\partial_1u\|_{L^2}\sup_{0\le \tau\le t}\|\Lambda_1^{-\sigma_1}u\|_{L^2}\\&\quad\,\,\cdot \left(\int_0^t(1+\tau)^{-\frac{\sigma_1\theta_1+\sigma_1+1}{s(2-\theta_1(1-\theta_2))}}\,\mathrm{d}\tau\right)^{\frac{2-\theta_1(1-\theta_2)}{2}},
\end{split}
\end{equation}
where the interpolation parameters are given by
\begin{equation}\label{11}
\theta_1=\frac{k+\sigma_1-1}{k+\sigma_1-\frac12},\qquad \theta_2=\frac{s+\sigma_1-\frac12}{s}.
\end{equation}
The assumptions on $k$, together with $s<2\sigma_1$ can guarantee 
$$
\theta_1>6\sigma_1-2>4s-2\sigma_1-2,
$$
which implies, 
$$
\frac{\sigma_1\theta_1+\sigma_1+1}
{s(2-\theta_1(1-\theta_2))}
>1.
$$
Applying Young's inequality to \eqref{J_1} we finally get 
\begin{equation}\label{eq2}
\begin{split}
\int_0^{t}|J_1|\,\mathrm{d}\tau&\lesssim E^{\frac12}_0(t)E^{\frac{\theta_1}{2}}_1(t)E^{\frac{1-\theta_1}{2}}_2(t)E^{\frac12}_3(t)\left(\int_0^t(1+\tau)^{-\frac{\sigma_1\theta_1+\sigma_1+1}{s(2-\theta_1(1-\theta_2))}}\,\mathrm{d}\tau\right)^{\frac{2-\theta_1(1-\theta_2)}{2}}\\&
 \lesssim  E^{\frac32}_0(t)+E^{\frac32}_1(t)+E^{\frac32}_2(t)+E^{\frac32}_3(t).
\end{split}
\end{equation}
{\bf{Estimate of $J_2$:}} By H\"older's inequality, Lemma \ref{2.2} and anisotropic type  interpolation inequalities, we deduce
\begin{equation}\label{J_22}
\begin{split}
|J_2|&\lesssim
 \|u_2\|_{L^{\frac{1}{\sigma_1}}_{x_1}L^{\infty}_{x_2}}\|\partial_2u\|_{L^2}\|\Lambda_1^{-\sigma_1}u\|_{L^2}\\&\lesssim \|\Lambda^{\frac12-\sigma_1}_1\Lambda^{\sigma}u_2\|^{\theta_3}_{L^2}\|\Lambda^{\frac12-\sigma_1}_1\Lambda^{k-\frac12+\sigma_1}u_2\|^{1-\theta_3}_{L^2}\|\partial_2u\|_{L^2}\|\Lambda_1^{-\sigma_1}u\|_{L^2}
\\&\lesssim \|\Lambda^{\sigma}u_2\|^{\theta_3\theta_4}_{L^2}\|\Lambda^s_1\Lambda^{\sigma}u_2\|^{\theta_3(1-\theta_4)}_{L^2}\|u_2\|^{1-\theta_3}_{H^k}\|u\|^{\theta_5}_{L^2}\|\partial^k_2u\|^{1-\theta_5}_{L^2}\|\Lambda_1^{-\sigma_1}u\|_{L^2}.
\end{split}
\end{equation}
Combining \eqref{J_22} with H\"older's inequality in time, we find
\begin{equation}\label{JJ2}
\begin{split}
&\int_0^{t}|J_2|\,\mathrm{d}\tau \\& \lesssim \sup_{0\le \tau\le t}(1+\tau)^{\frac{\theta_3\theta_4(\sigma_1+1)}{2s}} \|\Lambda^{\sigma}u_2\|^{\theta_3\theta_4}_{L^2} \left(\int_0^t(1+\tau)^{\frac{\sigma_1+1}{s}}\|\Lambda^s_1\Lambda^{\sigma}u_2\|^2_{L^2}\,\mathrm{d}\tau\right)^{\frac{\theta_3(1-\theta_4)}{2}}\\&\,\,\quad\cdot
\sup_{0\le \tau\le t}(1+\tau)^{\frac{\sigma_1\theta_5}{2s}}\|u\|^{\theta_5}_{L^2}\sup_{0\le \tau\le t}\|u\|^{2-\theta_3-\theta_5}_{H^k} \sup_{0\le \tau\le t}\|\Lambda_1^{-\sigma_1}u\|_{L^2}\\&\quad\,\,\cdot\left(\int_0^t(1+\tau)^{-\frac{(\sigma_1+1)\theta_3+\sigma_1\theta_5}{s(2-\theta_3(1-\theta_4))}}\,\mathrm{d}\tau\right)^{\frac{2-\theta_3(1-\theta_4)}{2}},
\end{split}
\end{equation}
where the interpolation parameters are given by
\begin{equation}\label{2222}
\theta_3=\frac{k+\sigma_1-1}{k+\sigma_1-\sigma-\frac12},\qquad \theta_4=\frac{s-\frac12+\sigma_1}{s},\qquad \theta_5=\frac{k-1}{k}.
\end{equation}
The assumptions on $k$ imply 
$$
\big(\frac32-3\sigma_1\big)k^2
+\big(-3\sigma_1^2+3\sigma\sigma_1+2\sigma_1-\frac32\big)k
+\sigma_1\big(\sigma+\frac12-\sigma_1\big)>0.
$$
Notice \eqref{2222} and $s<2\sigma_1$ we can guarantee 
$$
\frac32\theta_3+\sigma_1\theta_5>4\sigma_1>2s,
$$
which  yields
$$
\frac{(\sigma_1+1)\theta_3+\sigma_1\theta_5}
{s(2-\theta_3(1-\theta_4))}
>1.
$$
Therefore, the time integral part in \eqref{JJ2} is uniformly bounded in $t$. 
Using Young's inequality it's easy to get 
\begin{equation}\label{eq3}
\begin{split}
\int_0^{t}|J_2|\,\mathrm{d}\tau&\lesssim E^{\frac12}_0(t)E^{\frac{\theta_5}{2}}_1(t)E^{\frac{2-\theta_3-\theta_5}{2}}_2(t)E^{\frac{\theta_3}{2}}_5(t)\left(\int_0^t(1+\tau)^{-\frac{(\sigma_1+1)\theta_3+\sigma_1\theta_5}{s(2-\theta_3(1-\theta_4))}}\,\mathrm{d}\tau\right)^{\frac{2-\theta_3(1-\theta_4)}{2}}\\&
\lesssim
E_0^{\frac32}(t)
+E_1^{\frac32}(t)
+E_2^{\frac32}(t)
+E_5^{\frac32}(t).
\end{split}
\end{equation}
Combining \eqref{eq2} and \eqref{eq3} together, we then complete the proof of Lemma \ref{lem1}.
\end{proof}

Next, we turn to handle the higher-order energy estimate.  
The following lemma provides the bound control for $H^k$ norm of our solution.
\begin{lemma}\label{lem2}
Assume that 
$$k>\max\left\{\frac{\sigma_1}{3-6\sigma_1},\frac{12\sigma_1^2-13\sigma_1+4}{3-6\sigma_1},2\right\},\qquad k\in\mathbb N,$$
 and the energies are defined as in \eqref{eq1}, then we have
$$E_2(t)\lesssim E_2(0)+E^{\frac32}_1(t)+E^{\frac32}_2(t)+E^{\frac32}_3(t).$$
\begin{proof}
By the norm equivalence
$$
\|f\|^2_{H^k}\sim \|f\|^2_{L^2}+\sum^{2}_{m=1}\|\partial^k_mf\|^2_{L^2},
$$
it is sufficient to estimate the $L^2$ norm and the highest--order derivative norms. By a standard energy argument, we obtain the $L^2$ energy estimate for $u$,
$$
\|u(t)\|^2_{L^2}+2\nu\int_0^t\|\Lambda^s_1u(\tau)\|^2_{L^2}\,\mathrm{d}\tau=\|u_0\|^2_{L^2}.
$$
Applying $\partial_m^k$ to the equation \eqref{INS} and taking $L^2$ inner product with $\partial_m^ku$ lead to 
\begin{equation*}
\begin{split}
\frac12\sum^2_{m=1}\frac{\mathrm{d}}{\mathrm{d}t}\|\partial^k_mu\|^2_{L^2}+\nu\sum^2_{m=1}\|\Lambda^s_1\partial^k_mu\|^2_{L^2}&=-\sum^2_{m=1}<\partial^k_m(u\cdot\nabla u),\partial^k_mu>\\&\triangleq N_1+N_2+N_3,
\end{split}
\end{equation*}
where
\begin{equation*}
\begin{aligned}
N_1&=-\sum^2_{i,j=1}\sum^k_{\alpha=1}\binom{k}{\alpha}\int_{\mathbb{R}^2}\partial^{\alpha}_1u_j\partial_j\partial^{k-\alpha}_1u_i\partial^k_1u_i \,\mathrm{d}x,\\
N_2&=-\sum^2_{i=1}\sum^k_{\alpha=1}\binom{k}{\alpha}\int_{\mathbb{R}^2}\partial^{\alpha}_2u_1\partial_1\partial^{k-\alpha}_2u_i\partial^k_2u_i \,\mathrm{d}x,\\
N_3&=-\sum^2_{i=1}\sum^k_{\alpha=1}\binom{k}{\alpha}\int_{\mathbb{R}^2}\partial^{\alpha}_2u_2\partial_2\partial^{k-\alpha}_2u_i\partial^k_2u_i \,\mathrm{d}x\\&=\sum^2_{i=1}\sum^k_{\alpha=1}\binom{k}{\alpha}\int_{\mathbb{R}^2}\partial_1\partial^{\alpha-1}_2u_1\partial_2\partial^{k-\alpha}_2u_i\partial^k_2u_i \,\mathrm{d}x.
\end{aligned}
\end{equation*}
{\bf{Estimate of $N_1$:}}
Since this term contains enough derivatives in the $x_1$ direction, it can be estimated directly. For $1\leq \alpha\leq k-1$, H\"older's inequality and the Sobolev embedding theorem yield
\begin{equation*}
\begin{split}
&\left|
	\sum_{i,j=1}^{2}\sum_{\alpha=1}^{k-1}
	\binom{k}{\alpha}
	\int_{\mathbb{R}^2}
	\partial_1^\alpha u_j
	\partial_j\partial_1^{k-\alpha}u_i
	\partial_1^k u_i\,\mathrm{d}x
	\right|\\
 &\lesssim\sum^{2}_{i,j=1}\sum^{k-1}_{\alpha=1}\|\partial^{\alpha}_1u_j\|_{L^2_{x_1}L^{\infty}_{x_2}}\|\partial_j\partial^{k-\alpha}_1u_i\|_{L^4_{x_1}L^2_{x_2}}\|\partial^k_1u_i\|_{L^4_{x_1}L^2_{x_2}}\\&\lesssim  \sum^{2}_{i,j=1}\sum^{k-1}_{\alpha=1} \|\partial^{\alpha}_1u_j\|_{H^1}\|\Lambda^{\frac14}_1\partial_j\partial^{k-\alpha}_1u_i\|_{L^2}\|\Lambda^{\frac14}_1\partial^k_1u_i\|_{L^2}\\&\lesssim \sum^{2}_{i,j=1}\sum^{k-1}_{\alpha=1}\|\partial^{\alpha}_1u_j\|_{H^1}\|\Lambda^s_1\Lambda^{\frac14-s}_1\partial_j\partial^{k-\alpha}_1u_i\|_{L^2}\|\Lambda^s_1\Lambda^{\frac14-s}_1\partial^k_1u_i\|_{L^2}\\&\lesssim \|u\|_{H^k}\|\Lambda^s_1u\|^2_{H^k}.
\end{split}
\end{equation*}
The endpoint case $\alpha=k$ can be handled similarly like follows, 
\begin{equation*}
\begin{split}
\left|\sum^{2}_{i,j=1}\int_{\mathbb{R}^2} \partial^k_1u_j\partial_ju_i\partial^k_1u_i \,\mathrm{d}x\right|& \lesssim\sum^{2}_{i,j=1}\|\partial^k_1u_j\|_{L^4_{x_1}L^2_{x_2}}\|\partial_ju_i\|_{L^{2}_{x_1}L^{\infty}_{x_2}}\|\partial^k_1u_i\|_{L^4_{x_1}L^2_{x_2}}\\&\lesssim\|u\|_{H^k}\|\Lambda^s_1u\|^2_{H^k}.
\end{split}
\end{equation*}
Therefore,
\begin{equation}\label{3.21}
\int_0^{t}|N_1|\,\mathrm{d}\tau \lesssim \sup_{0\le \tau\le t} \|u\|_{H^k}\int_0^{t}\|\Lambda^s_1u\|^2_{H^k}\,\mathrm{d}\tau\lesssim  E^{\frac32}_2(t).
\end{equation}
{\bf{Estimates of $N_2$ and $N_3$:}} Despite the different distributions of derivatives, the terms $N_2$ and $N_3$ share the similar structure. 
They can be treated at the same time by distinguishing the cases of lower-order derivative terms and higher-order terms. 
Accordingly, we decompose them into {($l,m,n \in \{1,2\}$ )} two parts, 
\begin{equation*}
\begin{split}
N_{\mathrm{low}}^{l,m,n}&=\sum^{k-1}_{\beta=1}\int_{\mathbb{R}^2}\partial^{\beta}_2u_l\partial_1\partial^{k-\beta}_2u_m\partial^k_2u_n \,\mathrm{d}x,\\
N_{\mathrm{high}}^{l,m,n}&=\int_{\mathbb{R}^2}\partial^k_2u_l\partial_1u_m\partial^k_2u_n\,\mathrm{d}x.
\end{split}
\end{equation*}
We shall first deal with  $N_{\mathrm{low}}^{l,m,n}$. H\"older's inequality, Lemma \ref{2.3} and anisotropic interpolation inequalities give 
\begin{equation*}
\begin{split}
|N_{\mathrm{low}}^{l,m,n}|&\lesssim\sum^{k-1}_{\beta=1}\|\Lambda^s_1\partial^{k-\beta}_2u_m\|_{L^2}\|\Lambda^{1-s}_1(\partial^{\beta}_2u_l\partial^k_2u_n)\|_{L^2}\\&\lesssim \sum^{k-1}_{\beta=1}\|\Lambda^s_1u\|^{\widetilde{\theta}_1}_{L^2}\|\Lambda^s_1\partial^k_2u\|^{1-\widetilde{\theta}_1}_{L^2}\|\Lambda^{\frac34-\frac s2}_1\partial^{\beta}_2u\|_{L^2_{x_1}L^{\infty}_{x_2}}\|\Lambda^{\frac34-\frac s2}_1\partial^k_2u\|_{L^2}\\&\lesssim \sum^{k-1}_{\beta=1} \|u\|^{\widetilde{\theta}_1\widetilde{\theta}_2}_{L^2}\|\partial_1u\|^{\widetilde{\theta}_1(1-\widetilde{\theta}_2)}_{L^2}\|\Lambda^s_1\partial^k_2u\|^{1-\widetilde{\theta}_1}_{L^2}\|\Lambda^{\frac34-\frac s2}_1u\|^{\widetilde{\theta}_3}_{L^2}\|\Lambda^{\frac34-\frac s2}_1\partial^k_2u\|^{2-\widetilde{\theta}_3}_{L^2}\\&\lesssim \sum^{k-1}_{\beta=1} \|u\|^{\widetilde{\theta}_1\widetilde{\theta}_2}_{L^2}\|\partial_1u\|^{\widetilde{\theta}_1(1-\widetilde{\theta}_2)}_{L^2}\|\Lambda^s_1\partial^k_2u\|^{1-\widetilde{\theta}_1}_{L^2}\|u\|^{\widetilde{\theta}_3\widetilde{\theta}_4}_{L^2}\|\partial_1u\|^{\widetilde{\theta}_3(1-\widetilde{\theta}_4)}_{L^2}\\&\quad\,\,\cdot \|\partial^k_2u\|^{\widetilde{\theta}_5(2-\widetilde{\theta}_3)}_{L^2}\|\Lambda^s_1\partial^k_2u\|^{(1-\widetilde{\theta}_5)(2-\widetilde{\theta}_3)}_{L^2}\\&\lesssim\sum^{k-1}_{\beta=1} \|u\|^{\widetilde{\theta}_1\widetilde{\theta}_2+\widetilde{\theta}_3\widetilde{\theta}_4}_{L^2}\|\partial_1u\|^{\widetilde{\theta}_1(1-\widetilde{\theta}_2)+\widetilde{\theta}_3(1-\widetilde{\theta}_4)}_{L^2} \|\partial^k_2u\|^{\widetilde{\theta}_5(2-\widetilde{\theta}_3)}_{L^2}\|\Lambda^s_1\partial^k_2u\|^{1-\widetilde{\theta}_1+(1-\widetilde{\theta}_5)(2-\widetilde{\theta}_3)}_{L^2}\\&\lesssim \sum^{k-1}_{\beta=1} \|u\|^{\gamma_1}_{L^2}\|\partial_1u\|^{\gamma_2}_{L^2} \|\partial^k_2u\|^{\gamma_3}_{L^2}\|\Lambda^s_1\partial^k_2u\|^{\gamma_4}_{L^2}.
\end{split}
\end{equation*}
In which, 
$$
\widetilde\theta_1=\frac{\beta}{k},\qquad
\widetilde\theta_2=1-s,\qquad
\widetilde\theta_3=\frac{k-\frac12-\beta}{k},
$$
$$
\widetilde\theta_4=\frac{s}{2}+\frac14,\qquad
\widetilde\theta_5=\frac{\frac32s-\frac34}{s}, 
$$
and
$$
\gamma_1=\widetilde{\theta}_1\widetilde{\theta}_2+\widetilde{\theta}_3\widetilde{\theta}_4=\frac{\beta(1-s)}{k}+\frac{k-\frac12-\beta}{k}\left(\frac s2+\frac14\right),
$$
$$
\gamma_2=\widetilde{\theta}_1(1-\widetilde{\theta}_2)+\widetilde{\theta}_3(1-\widetilde{\theta}_4)=\frac{\beta s}{k}+\frac{k-\frac12-\beta}{k}\left(\frac34-\frac s2\right),
$$
$$
\gamma_3=\widetilde{\theta}_5(2-\widetilde{\theta}_3)=\left(\frac32-\frac{3}{4s}\right)\left(1+\frac{\beta+\frac12}{k}\right),
$$
$$
\gamma_4=1-\widetilde{\theta}_1+(1-\widetilde{\theta}_5)(2-\widetilde{\theta}_3)=1-\frac{\beta}{k}+\left(\frac{3}{4s}-\frac12\right)\left(1+\frac{\beta+\frac12}{k}\right).
$$
Moreover, the assumptions on $k$ together with $s<2\sigma_1$ imply that 
$$
k>\frac{\sigma_1}{2(\sigma_1+\frac32-2s)},
$$
which ensures
$$
\frac{\sigma_1\gamma_1+(\sigma_1+1)\gamma_2}{s(2-\gamma_4)}>1.
$$
A direct computation gives $\gamma_1+\gamma_2+\gamma_3+\gamma_4=3$.
Since $0<\gamma_4<2$, applying H\"older's inequality in time then yields 
\begin{equation}\label{3.22}
\begin{aligned}
\int_0^t |N_{\mathrm{low}}^{l,m,n}|\,\mathrm{d}\tau
&\lesssim
\sum_{\beta=1}^{k-1}
\sup_{0\leq\tau\leq t}
\left(
(1+\tau)^{
\frac{\sigma_1\gamma_1+(\sigma_1+1)\gamma_2}{2s}}
\|u\|_{L^2}^{\gamma_1}
\|\partial_1u\|_{L^2}^{\gamma_2}
\|\partial_2^ku\|_{L^2}^{\gamma_3}
\right)
\\
&\qquad\cdot
\left(
\int_0^t
\|\Lambda_1^s\partial_2^ku\|_{L^2}^2
\,\mathrm{d}\tau
\right)^{\frac{\gamma_4}{2}}
\left(
\int_0^t
(1+\tau)^{
-\frac{\sigma_1\gamma_1+(\sigma_1+1)\gamma_2}
{s(2-\gamma_4)}}
\,\mathrm{d}\tau
\right)^{\frac{2-\gamma_4}{2}}
\\
&\lesssim
\sum_{\beta=1}^{k-1}
E_1^{\frac{\gamma_1}{2}}(t)
E_3^{\frac{\gamma_2}{2}}(t)
E_2^{\frac{\gamma_3+\gamma_4}{2}}(t)
\\
&\lesssim
E_1^{\frac32}(t)
+
E_2^{\frac32}(t)
+
E_3^{\frac32}(t).
\end{aligned}
\end{equation}
For the endpoint case $\beta=k$, we denote the corresponding term by $N_{\mathrm{high}}^{l,m,n}$.
Using H\"older's inequality, interpolation inequality and Young's inequality, we shall derive
\begin{equation}\label{wang}
\begin{split}
|N_{\mathrm{high}}^{l,m,n}|
&\lesssim
\|\partial_1u\|_{L^2_{x_1}L^\infty_{x_2}}
\|\partial_2^ku\|_{L^4_{x_1}L^2_{x_2}}^2
\\
&\lesssim
\|\partial_1u\|_{L^2}^{\frac12}
\|\partial_1\partial_2u\|_{L^2}^{\frac12}
\|\Lambda_1^{\frac14}\partial_2^ku\|_{L^2}^2
\\
&\lesssim
\|\partial_1u\|_{L^2}^{\frac12}
\|\partial_1u\|_{L^2}^{
 \frac{k+s-2}{2(k+s-1)}}
\|\Lambda_2^{k+s-1}\partial_1u\|_{L^2}^{
 \frac{1}{2(k+s-1)}}
\|\partial_2^ku\|_{L^2}^{2-\frac{1}{2s}}
\|\Lambda_1^s\partial_2^ku\|_{L^2}^{\frac{1}{2s}}
\\
&\lesssim
\|\partial_1u\|_{L^2}^{
 \frac{2k+2s-3}{2(k+s-1)}}
\|u\|_{H^k}^{
 2-\frac{3}{4s}-\frac{1}{4(k+s-1)}}
\|u\|_{H^k}^{
 \frac{k+2s-1}{4s(k+s-1)}}
\|\Lambda_1^su\|_{H^k}^{
 \frac{k+2s-1}{2s(k+s-1)}}
\\
&\lesssim
\left(
\|\partial_1u\|_{L^2}^{
 \frac{2k+2s-3}{2(k+s-1)}}
\|u\|_{H^k}^{
 2-\frac{3}{4s}-\frac{1}{4(k+s-1)}}
\right)^{
 \frac{4s(k+s-1)}
 {4s(k+s-1)-k-2s+1}}
+
\|u\|_{H^k}\|\Lambda_1^su\|_{H^k}^2.
\end{split}
\end{equation}
The assumptions on $k$ together with $s<2\sigma_1$ will imply 
$$
k>
\frac{4s^2-2\sigma_1s-8s+3\sigma_1+4}
{2\sigma_1+3-4s}.
$$
Consequently,
\begin{equation}\label{3.23}
\begin{aligned}
\int_0^t |N_{high}^{l,m,n}|\,\mathrm{d}\tau
&\lesssim
\sup_{0\leq\tau\leq t}(1+\tau)^{\frac{(\sigma_1+1)(2k+2s-3)}{4s(k+s-1)-k-2s+1}}\|\partial_1u\|_{L^2}^{\frac{2s(2k+2s-3)}{4s(k+s-1)-k-2s+1}}\\&\,\,\quad\cdot
\sup_{0\le\tau\le t}
\|u(\tau)\|_{H^k}^{(2-\frac{3}{4s}-\frac{1}{4(k+s-1)})(\frac{4s(k+s-1)}
 {4s(k+s-1)-k-2s+1})}
 \\&\,\,\quad\cdot\int_0^t (1+\tau)^{-\frac{(\sigma_1+1)(2k+2s-3)}{4s(k+s-1)-k-2s+1}}\,\mathrm{d}\tau+\int_0^t\|u\|_{H^k}\|\Lambda^s_1u\|^2_{H^k}\,\mathrm{d}\tau
\\&\lesssim
E^{\frac32}_2(t)+E^{\frac32}_3(t),
\end{aligned}
\end{equation}
It follows from \eqref{3.21}, \eqref{3.22} and \eqref{3.23} that the proof of Lemma \ref{lem2} finishes.
\end{proof}

\end{lemma}

Based on the previous energy estimate, we now turn to deal with the time-weighted energy  for $\partial_1u$.
\begin{lemma}\label{lem3}
Assume that
$$k>\max{\left\{\frac{\frac{11}{4}+7\sigma_1-4\sigma+\sigma_1(2\sigma_1+\sigma)}{1-2\sigma_1},\frac{12\sigma_1^2-13\sigma_1+4}{3-6\sigma_1},2\right\}},$$
and the energies are defined as in \eqref{eq1}, then we have
$$
E_3(t)\lesssim E_3(0)+E^{\frac{s}{1+\sigma_1}}_0(t)E^{\frac{1+\sigma_1-s}{1+\sigma_1}}_3(t)+ E^{\frac32}_1(t)+E^{\frac32}_2(t)+E^{\frac32}_3(t)+E^{\frac32}_5(t).
$$
\begin{proof}
Applying $\partial_1$ to equation \eqref{INS}, taking the $L^2$ inner product with $\partial_1u$, and multiplying by the time weight $(1+t)^{\frac{\sigma_1+1}{s}}$, we shall get 
\begin{equation*}
\frac12\frac{\mathrm{d}}{\mathrm{d}t}\left((1+t)^{\frac{\sigma_1+1}{s}}\|\partial_1u\|^2_{L^2}\right)+\nu(1+t)^{\frac{\sigma_1+1}{s}}\|\Lambda^{1+s}_1u\|^2_{L^2}\triangleq I_1+I_2,
\end{equation*}
where
\begin{equation*}
\begin{aligned}
I_1&=\frac12\frac{\sigma_1+1}{s}(1+t)^{\frac{\sigma_1+1-s}{s}}\|\partial_1u\|^2_{L^2},\\
I_2&=-(1+t)^{\frac{\sigma_1+1}{s}}<\partial_1(u\cdot\nabla u),\partial_1u>\\&=-(1+t)^{\frac{\sigma_1+1}{s}}<\partial_1u_1\partial_1u,\partial_1u>-(1+t)^{\frac{\sigma_1+1}{s}}<\partial_1u_2\partial_2u,\partial_1u>.
\end{aligned}
\end{equation*}
{\bf{Estimate of $I_1$:}} Interpolation inequality and H\"older's inequality in time give
\begin{equation}\label{I_1}
\begin{split}
\int_0^t|I_1|\,\mathrm{d}\tau& \lesssim\int_0^t(1+\tau)^{\frac{\sigma_1+1-s}{s}}\|\Lambda^{s-\sigma_1}_1u\|^{\frac{2s}{1+\sigma_1}}_{L^2}\|\Lambda^{1+s}_1u\|_{L^2}^{\frac{2(1+\sigma_1-s)}{1+\sigma_1}}\,\mathrm{d}\tau\\&\lesssim \int_0^t\|\Lambda^{s-\sigma_1}_1u\|^{\frac{2s}{1+\sigma_1}}_{L^2}\left((1+\tau)^{\frac{\sigma_1+1}{2s}}\|\Lambda^{1+s}_1u\|_{L^2}\right)^{\frac{2(1+\sigma_1-s)}{1+\sigma_1}}\,\mathrm{d}\tau\\&\lesssim E^{\frac{s}{1+\sigma_1}}_0(t)E^{\frac{1+\sigma_1-s}{1+\sigma_1}}_3(t).
\end{split}
\end{equation}
{\bf{Estimate of $I_2$:}} Using H\"older's inequality, it's easy to get 
\begin{equation}\label{22}
\begin{split}
\int_0^t|I_2|\,\mathrm{d}\tau& \lesssim \int_0^t(1+\tau)^{\frac{\sigma_1+1}{s}}\|\partial_1u_1\|_{L^2_{x_1}L^{\infty}_{x_2}}
\|\partial_1u\|_{L^4_{x_1}L^2_{x_2}}^2\,\mathrm{d}\tau\\&\,\,\quad+\sup_{0\le \tau\le t}(1+\tau)^{\frac{\sigma_1+1}{2s}}\|\partial_1u\|_{L^2}\int_0^t(1+\tau)^{\frac{\sigma_1+1}{2s}}\|\partial_1u_2\partial_2u\|_{L^2}\,\mathrm{d}\tau\\&\triangleq I_{21}+I_{22}.
\end{split}
\end{equation}
Using the same anisotropic type 
interpolation like in \eqref{wang} will lead to 
\begin{equation*}
\begin{aligned}
&(1+\tau)^{\frac{1+\sigma_1}{s}}
\|\partial_1u_1\|_{L^2_{x_1}L^\infty_{x_2}}
\|\partial_1u\|_{L^4_{x_1}L^2_{x_2}}^2
\\&\lesssim
(1+\tau)^{\frac{1+\sigma_1}{s}}
\|\partial_1u\|_{L^2}^{
 1-\frac1{2(k+s-1)}}
\|\Lambda_1^su\|_{H^k}^{
 \frac1{2(k+s-1)}}
\|\partial_1u\|_{L^2}^{2-\frac1{2s}}
\|\Lambda_1^{1+s}u\|_{L^2}^{\frac1{2s}}\\
&\lesssim
 (1+\tau)^{\frac{1+\sigma_1}{s}}
\|\partial_1u\|_{L^2}^{
 3-\frac1{2(k+s-1)}-\frac1{2s}}
\|\Lambda_1^su\|_{H^k}^{
 \frac1{2(k+s-1)}}
\|\Lambda_1^{1+s}u\|_{L^2}^{\frac1{2s}}.
\end{aligned}
\end{equation*}
The assumptions on $k$, together with $s<2\sigma_1$, 
	ensure that the last time integral below is uniformly
	bounded in $t$.
Thus, H\"older's inequality in time yields
\begin{equation}\label{I_21}
\begin{split}
I_{21}
&\lesssim
\int_0^t
(1+\tau)^{\frac{1+\sigma_1}{s}}
\|\partial_1u\|_{L^2}^{
3-\frac{1}{2(k+s-1)}-\frac{1}{2s}}
\|\Lambda_1^su\|_{H^k}^{\frac{1}{2(k+s-1)}}
\|\Lambda_1^{1+s}u\|_{L^2}^{\frac{1}{2s}}
\,\mathrm{d}\tau
\\
&\lesssim
\sup_{0\leq\tau\leq t}
\left(
(1+\tau)^{\frac{1+\sigma_1}{s}}
\|\partial_1u\|_{L^2}^{2}
\right)^{
\frac32-\frac{1}{4(k+s-1)}-\frac{1}{4s}}
\left(
\int_0^t
\|\Lambda_1^su\|_{H^k}^{2}\,\mathrm{d}\tau
\right)^{\frac{1}{4(k+s-1)}}
\\
&\quad\,\cdot
\left(
\int_0^t
(1+\tau)^{\frac{1+\sigma_1}{s}}
\|\Lambda_1^{1+s}u\|_{L^2}^{2}\,\mathrm{d}\tau
\right)^{\frac{1}{4s}}\\
&\quad\,\cdot
\left(
\int_0^t
(1+\tau)^{
-\frac{(1+\sigma_1)(2k+2s-3)}
{4s(k+s-1)-k-2s+1}}
\,\mathrm{d}\tau
\right)^{
1-\frac{1}{4(k+s-1)}-\frac{1}{4s}}
\\
&\lesssim
E_3^{
	\frac32-\frac{1}{4(k+s-1)}}(t)
E_2^{\frac{1}{4(k+s-1)}}(t)
\left(
\int_0^t
(1+\tau)^{
-\frac{(1+\sigma_1)(2k+2s-3)}
{4s(k+s-1)-k-2s+1}}
\,\mathrm{d}\tau
\right)^{
1-\frac{1}{4(k+s-1)}-\frac{1}{4s}}\\&\lesssim E_{2}^{\frac32}(t)+E_{3}^{\frac32}(t).
\end{split}
\end{equation}
For the term $I_{22}$,  we first deal with the part 
	$\|\partial_1u_2\partial_2u\|_{L^2}$. By  H\"older's inequality and anisotropic interpolation there is, 
\begin{equation*}
\begin{split}
\|\partial_1u_2\partial_2u\|_{L^2}&\lesssim \|\partial_1u_2\|_{L^2_{x_1}L^{\infty}_{x_2}}\|\partial_2u\|_{L^{\infty}_{x_1}L^2_{x_2}}\\&\lesssim \|\partial_1\Lambda^{\sigma}_2u_2\|^{\bar{\theta}_1}_{L^2}\|\partial_1\Lambda^{k-1}_2u_2\|^{1-\bar{\theta}_1}_{L^2}\|\partial_2u\|^{\frac12}_{L^2}\|\partial_1\partial_2u\|^{\frac12}_{L^2}\\&\lesssim \|\Lambda^s_1\Lambda^{\sigma}u_2\|^{\bar{\theta}_1\bar{\theta}_2}_{L^2}\|\Lambda^{k-\sigma}_1\Lambda^{\sigma}_2u_2\|^{\bar{\theta}_1(1-\bar{\theta}_2)}_{L^2}\|u\|^{1-\bar{\theta}_1}_{H^k}\|u\|^{\frac{\bar{\theta}_3}{2}}_{L^2}\|u\|^{\frac{1-\bar{\theta}_3}{2}}_{H^k}\\&\,\,\quad\cdot \|\partial_1u\|^{\frac{\bar{\theta}_4}{2}}_{L^2}\|u\|^{\frac{1-\bar{\theta}_4}{2}}_{H^k}\\&\lesssim \|\Lambda^s_1\Lambda^{\sigma}u_2\|^{\bar{\theta}_1\bar{\theta}_2}_{L^2}\|u\|^{\frac{\bar{\theta}_3}{2}}_{L^2}\|\partial_1u\|^{\frac{\bar{\theta}_4}{2}}_{L^2}\|u\|^{\bar{\theta}_1(1-\bar{\theta}_2)+1-\bar{\theta}_1+\frac{1-\bar{\theta}_3}{2}+\frac{1-\bar{\theta}_4}{2}}_{H^k},
\end{split}
\end{equation*}
where the interpolation exponents are determined by
$$
\bar{\theta}_1=\frac{k-\frac32}{k-1-\sigma},\qquad\bar{\theta}_2=\frac{k-\sigma-1}{k-\sigma-s},
$$
$$
 \bar{\theta}_3=\frac{k-1}{k},\qquad \bar{\theta}_4=\frac{k-2}{k-1}.
$$
Furthermore, applying Hölder's inequality we shall derive 
\begin{equation}\label{ddd}
\begin{split}
&\int_0^t(1+\tau)^{\frac{\sigma_1+1}{2s}}\|\partial_1u_2\partial_2u\|_{L^2}\,\mathrm{d}\tau\\&\lesssim \sup_{0\le \tau\le t}\|u\|^{\bar{\theta}_1(1-\bar{\theta}_2)+1-\bar{\theta}_1+\frac{1-\bar{\theta}_3}{2}+\frac{1-\bar{\theta}_4}{2}}_{H^k}\sup_{0\le \tau\le t}(1+\tau)^{\frac{\sigma_1\bar{\theta}_3}{4s}}\|u\|^{\frac{\bar{\theta}_3}{2}}_{L^2}\\&\,\,\quad\cdot \sup_{0\le \tau\le t}(1+\tau)^{\frac{\bar{\theta}_4(\sigma_1+1)}{4s}}\|\partial_1u\|^{\frac{\bar{\theta}_4}{2}}_{L^2}\left(\int_0^t(1+\tau)^{\frac{\sigma_1+1}{s}}\|\Lambda^s_1\Lambda^{\sigma}u_2\|^2_{L^2}\,\mathrm{d}\tau\right)^{\frac{\bar{\theta}_1\bar{\theta}_2}{2}}\\&\,\,\quad\cdot\left(\int_0^t(1+\tau)^{-\frac{2(\sigma_1+1)\bar{\theta}_1\bar{\theta}_2+\sigma_1\bar{\theta}_3+(\sigma_1+1)\bar{\theta}_4-2(\sigma_1+1)}{2s(2-\bar{\theta}_1\bar{\theta}_2)}}\,\mathrm{d}\tau\right)^{\frac{2-\bar{\theta}_1\bar{\theta}_2}{2}}.
\end{split}
\end{equation}
To verify the integrability of the last time integral part in \eqref{ddd}, we set 
\begin{equation*}
\begin{aligned}
P(s,k)
&=(1+2\sigma_1-2s)k^3
+(4s^2+4s\sigma+\sigma-5-7\sigma_1)k^2
\\
&\quad
+(2\sigma_1s+2\sigma_1\sigma+4\sigma_1
-4s^2-4s\sigma+3s+3)k
-\sigma_1(s+\sigma).
\end{aligned}
\end{equation*}
Notice the assumptions on $k$ we have 
$$
P(s,k)>(1-2\sigma_1)k^3-\big(\frac{11}{4}+7\sigma_1-4\sigma+\sigma_1(2\sigma_1+\sigma)\big)k^2>0.
$$
Consequently,
$$
\frac{2(\sigma_1+1)\bar{\theta}_1\bar{\theta}_2+\sigma_1\bar{\theta}_3+(\sigma_1+1)\bar{\theta}_4-2(\sigma_1+1)}{2s(2-\bar{\theta}_1\bar{\theta}_2)}>1.
$$
Hence, the last time integral part in \eqref{ddd}
	is uniformly bounded in $t$. Moreover, simple calculation yields 
$$
\frac{\bar{\theta}_1\bar{\theta}_2}{2}+\frac{\bar{\theta}_3}{4}+\frac{\bar{\theta}_4}{4}+\frac{\bar{\theta}_1(1-\bar{\theta}_2)+1-\bar{\theta}_1+\frac{1-\bar{\theta}_3}{2}+\frac{1-\bar{\theta}_4}{2}}{2}=1.
$$
Applying Young's inequality to \eqref{ddd} we then obtain 
\begin{equation}\label{I_22}
\begin{split}
I_{22}&\lesssim \sup_{0\le \tau\le t}(1+\tau)^{\frac{\sigma_1+1}{2s}}\|\partial_1u\|_{L^2}\int_0^t(1+\tau)^{\frac{\sigma_1+1}{2s}}\|\partial_1u_2\partial_2u\|_{L^2}\,\mathrm{d}\tau.
\\&
\lesssim E^{\frac{\bar{\theta}_1\bar{\theta}_2}{2}}_5(t)E^{\frac{\bar{\theta}_3}{4}}_1(t)E^{\frac{\bar{\theta}_4+2}{4}}_3(t) 
E^{\frac{\bar{\theta}_1(1-\bar{\theta}_2)+1-\bar{\theta}_1+\frac{1-\bar{\theta}_3}{2}+\frac{1-\bar{\theta}_4}{2}}{2}}_2(t)\\&\,\,\quad\cdot\left(\int_0^t(1+\tau)^{-\frac{2(\sigma_1+1)\bar{\theta}_1\bar{\theta}_2+\sigma_1\bar{\theta}_3+(\sigma_1+1)\bar{\theta}_4-2(\sigma_1+1)}{2s(2-\bar{\theta}_1\bar{\theta}_2)}}\,\mathrm{d}\tau\right)^{\frac{2-\bar{\theta}_1\bar{\theta}_2}{2}}\\&\lesssim E^{\frac32}_1(t)+E^{\frac32}_2(t)+E^{\frac32}_3(t)+E^{\frac32}_5(t).
\end{split}
\end{equation}
Combining \eqref{I_1}, \eqref{I_21} and \eqref{I_22} together,  we conclude the proof of Lemma \ref{lem3}.
\end{proof}

\end{lemma}

Due to the anisotropic nature of the dissipation, the second velocity component
requires separated analysis. Although the dissipation acts only in the $x_1$
direction, the divergence-free structure provides additional property for $u_2$,
which enables us to obtain improved decay estimate. We shall begin with the negative-order derivatives for $u_2$.
\begin{lemma}\label{lem4}
Assume that
$$k>\max\left\{\frac{\sigma_1}{2(1-2\sigma_1)},
\frac{6\sigma_1^2-6\sigma_1+2}{3-6\sigma_1},
\frac{\left|-6\sigma_1^2+6\sigma\sigma_1+4\sigma_1-3\right|}
{3-6\sigma_1},\frac{\sigma_1}{3-4\sigma_1},
1
\right\},$$ 
and the energies are defined as in \eqref{eq1}, then we have
$$
E_4(t)\lesssim E_4(0)+E^{\frac32}_1(t)+E^{\frac32}_2(t)+E^{\frac32}_3(t)+E^{\frac32}_4(t)+E^{\frac32}_5(t).
$$
\end{lemma}
\begin{proof}
Recall the vorticity equation $(\omega=\nabla^{\perp}\cdot u)$ here 
$$
\partial_t\omega+\nu\Lambda^{2s}_1\omega+\nabla^{\perp} \cdot (u\cdot\nabla u)=0.
$$
Since $u$ is divergence-free, there exists a stream function $\phi$ such that $u=\nabla^{\perp}\phi$, and the Biot-Savart law yields {$u=-\nabla^{\perp}(-\Delta)^{-1}\omega.$} Consequently, the second velocity component $u_2$ satisfies the identity 
\begin{equation}\label{u_2}
\partial_tu_2+\nu\Lambda^{2s}_1u_2= (-\Delta)^{-1}\partial_1\nabla^{\bot}\cdot(u\cdot\nabla u).
\end{equation}
Applying the operator $\Lambda_1^{-\sigma_1-1}\Lambda^\sigma$ to
\eqref{u_2} and taking the $L^2$ inner product with
$\Lambda_1^{-\sigma_1-1}\Lambda^\sigma u_2$, we obtain
\begin{equation*}
\frac12\frac{\mathrm{d}}{\mathrm{d}t}\|\Lambda^{-\sigma_1-1}_1\Lambda^{\sigma}u_2\|^2_{L^2}+\nu\|\Lambda^{-\sigma_1-1+s}_1\Lambda^{\sigma}u_2\|^2_{L^2}= Q,
\end{equation*}
where
\begin{equation*}
\begin{aligned}
Q&=\langle\Lambda^{-\sigma_1-1}_1\Lambda^{\sigma}[(-\Delta)^{-1}\partial_1\nabla^{\bot}\cdot(u\cdot\nabla u)],\Lambda_1^{-\sigma_1-1}\Lambda^{\sigma}u_2\rangle,\\&=-<\Lambda^{-\sigma_1-1}_1\Lambda^{\sigma}[(-\Delta)^{-1}\partial_1\partial_2\partial_1(u_1u_1)],\Lambda_1^{-\sigma_1-1}\Lambda^{\sigma}u_2>\\&\,\,\quad-<\Lambda^{-\sigma_1-1}_1\Lambda^{\sigma}[(-\Delta)^{-1}\partial_1\partial_2\partial_2(u_2u_1)],\Lambda_1^{-\sigma_1-1}\Lambda^{\sigma}u_2>\\&\,\,\quad+<\Lambda^{-\sigma_1-1}_1\Lambda^{\sigma}[(-\Delta)^{-1}\partial_1\partial_1\partial_1(u_1u_2)],\Lambda_1^{-\sigma_1-1}\Lambda^{\sigma}u_2>\\&\,\,\quad+<\Lambda^{-\sigma_1-1}_1\Lambda^{\sigma}[(-\Delta)^{-1}\partial_1\partial_1\partial_2(u_2u_2)],\Lambda_1^{-\sigma_1-1}\Lambda^{\sigma}u_2>,\\&\triangleq Q_1+Q_2+Q_3+Q_4.
\end{aligned}
\end{equation*}
{\bf{Estimate of $Q_1$:}} By the $L^2$ boundness of Calderón--Zygmund operators, we can obtain 
\begin{equation*}
\begin{split}
|Q_1|&\lesssim\|\Lambda^{-\sigma_1-1}_1\Lambda^{\sigma}(-\Delta)^{-1}\partial_1\partial_2(\partial_1u_1u_1)\|_{L^2}\|\Lambda_1^{-\sigma_1-1}\Lambda^{\sigma}u_2\|_{L^2}\\&\lesssim\|\Lambda^{-\sigma_1}_1\Lambda^{\sigma-1}(\partial_1u_1u_1)\|_{L^2}\|\Lambda_1^{-\sigma_1-1}\Lambda^{\sigma}u_2\|_{L^2}. 
\end{split}
\end{equation*}
To redistribute the negative derivatives, we choose
\[
\frac12-\sigma<\delta<\frac12-\sigma_1.
\]
{This interval is nonempty since the assumption $\sigma>\sigma_1$. Put
$\alpha=\sigma_1+\delta$ and $\beta=1-\sigma-\delta$, then
$0<\alpha,\beta<\frac{1}{2}$. The pointwise multiplier comparison}
\[
|\xi_1|^{-\sigma_1}|\xi|^{\sigma-1}
\le |\xi_1|^{-\sigma_1-\delta}|\xi|^{\sigma-1+\delta}
\le |\xi_1|^{-\alpha}|\xi_2|^{-\beta}
\]
coming from {Lemma~\ref{2.1}
and one-dimensional Gagliardo--Nirenberg inequality}, yields that 
\[
\begin{aligned}
|Q_1|&\lesssim
\|\Lambda_1^{-\alpha}\Lambda_2^{-\beta}(\partial_1u_1u_1)\|_{L^2}
\|\Lambda_1^{-\sigma_1-1}\Lambda^\sigma u_2\|_{L^2}\\
&\lesssim\Bigl(
\|\partial_1u_1u_1\|_{L^1}+\|\partial_1u_1u_1\|_{L^2}\\
&\qquad+\|\partial_1u_1u_1\|_{L^1_{x_1}L^2_{x_2}}
+\|\partial_1u_1u_1\|_{L^1_{x_2}L^2_{x_1}}\Bigr)
\|\Lambda_1^{-\sigma_1-1}\Lambda^\sigma u_2\|_{L^2}\\
&\lesssim\Bigl(
\|\partial_1u_1\|_{L^2}\|u_1\|_{L^2}
+\|\partial_1u_1\|_{L^2}^{1-\frac1{2s}}
\|\Lambda_1^{1+s}u_1\|_{L^2}^{\frac1{2s}}
\|u_1\|_{L^2}^{1-\frac1{2k}}\|\partial_2^ku_1\|_{L^2}^{\frac1{2k}}\\
&\qquad+\|\partial_1u_1\|_{L^2}
\|u_1\|_{L^2}^{1-\frac1{2k}}\|\partial_2^ku_1\|_{L^2}^{\frac1{2k}}
+\|\partial_1u_1\|_{L^2}^{\frac32}\|u_1\|_{L^2}^{\frac12}\Bigr)
\|\Lambda_1^{-\sigma_1-1}\Lambda^\sigma u_2\|_{L^2}.
\end{aligned}
\]
Integrating with respect to time gives
\begin{equation*}
\begin{split}
\int_0^t|Q_1|\,\mathrm{d}\tau&\lesssim \sup_{0\le \tau\le t}(1+\tau)^{\frac{\sigma_1+1}{2s}} \|\partial_1u_1\|_{L^2}\sup_{0\le \tau\le t}(1+\tau)^{\frac{\sigma_1}{2s}}\|u_1\|_{L^2} \\
&\,\,\quad\cdot\sup_{0\le \tau\le t}\|\Lambda_1^{-\sigma_1-1}\Lambda^{\sigma}u_2\|_{L^2}\int_0^t (1+\tau)^{-\frac{2\sigma_1+1}{2s}}\,\mathrm{d}\tau\\&\,\,\quad+ \sup_{0\le \tau\le t}(1+\tau)^{\frac{(\sigma_1+1)(2s-1)}{4s^2}}\|\partial_1u_1\|^{1-\frac{1}{2s}}_{L^2}\left(\int_0^t(1+\tau)^{\frac{\sigma_1+1}{s}}\|\Lambda^{1+s}_1u_1\|_{L^2}^2\,\mathrm{d}\tau\right)^{\frac{1}{4s}}\\&\,\,\quad\cdot\sup_{0\le \tau\le t}(1+\tau)^{\frac{\sigma_1(2k-1)}{4sk}}\|u_1\|^{1-\frac{1}{2k}}_{L^2} \sup_{0\le \tau\le t}\|u\|_{H^k}^{\frac{1}{2k}}
\\&\,\,\quad\cdot\sup_{0\le \tau\le t}\|\Lambda_1^{-\sigma_1-1}\Lambda^{\sigma}u_2\|_{L^2}\left(\int_0^t(1+\tau)^{-\frac{2+4\sigma_1-\frac{\sigma_1}{k}}{4s-1}}\,\mathrm{d}\tau\right)^{\frac{4s-1}{4s}}\\&\,\,\quad+\sup_{0\le \tau\le t}(1+\tau)^{\frac{\sigma_1+1}{2s}}\|\partial_1u_1\|_{L^2}\sup_{0\le \tau\le t}(1+\tau)^{\frac{\sigma_1(k-\frac12)}{2sk}}\|u_1\|_{L^2}^{1-\frac{1}{2k}}\sup_{0\le \tau\le t}\|\partial^k_2u_1\|^{\frac{\frac12}{k}}_{L^2}\\&\,\,\quad\cdot\sup_{0\le \tau\le t}\|\Lambda_1^{-\sigma_1-1}\Lambda^{\sigma}u_2\|_{L^2} \int_0^t(1+\tau)^{-\frac{k(\sigma_1+1)+\sigma_1(k-\frac12)}{2sk}} \,\mathrm{d}\tau\\&\,\,\quad+ \sup_{0\le \tau\le t}(1+\tau)^{\frac{3\sigma_1+3}{4s}}\|\partial_1u_1\|^{\frac32}_{L^2}\sup_{0\le \tau\le t}(1+\tau)^{\frac{\sigma_1}{4s}}\|u_1\|^{\frac12}_{L^2}\\
&\,\,\quad\cdot\sup_{0\le \tau\le t}\|\Lambda_1^{-\sigma_1-1}\Lambda^{\sigma}u_2\|_{L^2}\int_0^t(1+\tau)^{-\frac{4\sigma_1+3}{4s}} \,\mathrm{d}\tau.
\end{split}
\end{equation*}
The assumptions on $k$ together with $\frac{3}{4}<s<2\sigma_1$ imply that 
$$
\frac{2\sigma_1+1}{2s}>1,\qquad \frac{2+4\sigma_1-\frac{\sigma_1}{k}}{4s-1}>1,
$$
$$
\frac{k(\sigma_1+1)+\sigma_1(k-\frac12)}{2sk}>1,\qquad\frac{4\sigma_1+3}{4s}>1.
$$
Consequently, we obtain 
\begin{equation}\label{4.1}
\begin{split}
\int_0^t|Q_1|\,\mathrm{d}\tau &\lesssim E^{\frac12}_3(t)E^{\frac12}_1(t)E^{\frac12}_4(t)+E^{\frac12}_3(t)E^{\frac12-\frac{1}{4k}}_1(t)E^{\frac{1}{4k}}_2(t)E^{\frac12}_4(t)\\
&\qquad\qquad+E^{\frac34}_3(t)E^{\frac14}_1(t)E^{\frac12}_4(t)\\& \lesssim E^{\frac32}_1(t)+E^{\frac32}_2(t)+E^{\frac32}_3(t)+E^{\frac32}_4(t).
\end{split}
\end{equation}
{\bf Estimate of $Q_2+Q_3$:}
Notice the $L^2$ boundness of $(-\Delta)^{-1}\partial_i\partial_j$ and
Lemma{~\ref{lem:negative-product-1}} with $a=\sigma_1$, $b=\sigma$,  there is 
\[
\begin{aligned}
|Q_2+Q_3|&\lesssim
\|\Lambda_1^{-\sigma_1}\Lambda^\sigma(u_1u_2)\|_{L^2}
\|\Lambda_1^{-\sigma_1-1}\Lambda^\sigma u_2\|_{L^2}\\
&\lesssim\Bigl(
\|\Lambda_1^{\frac12-\sigma_1}u_1^\sharp\|_{L^2_{x_1}L^\infty_{x_2}}
\|\Lambda^\sigma u_2^\sharp\|_{L^2}\\
&\qquad+
\|\Lambda_1^{\frac12-\sigma_1}u_2^\sharp\|_{L^2_{x_1}L^\infty_{x_2}}
\|\Lambda^\sigma u_1^\sharp\|_{L^2}\Bigr)
\|\Lambda_1^{-\sigma_1-1}\Lambda^\sigma u_2\|_{L^2}.
\end{aligned}
\]
Notice again $\|u\|_{L^2} = \|u^\sharp\|_{L^2}$, using the similar method as in \eqref{J_11} and \eqref{J_22}, we obtain 

\begin{equation*}
\begin{split}
\|\Lambda^{\frac12-\sigma_1}_1 u_1^\sharp\|_{L^2_{x_1}L^{\infty}_{x_2}}\|\Lambda^{\sigma}u_2^\sharp\|_{L^2}&\lesssim \|\Lambda^{\frac12-\sigma_1}_1u_1\|^{\theta_1}_{L^2}\|\Lambda^{\frac12-\sigma_1}_1\Lambda^{k-\frac12+\sigma_1}_2u_1\|^{1-\theta_1}_{L^2}\|\Lambda^{\sigma}u_2\|_{L^2}\\&\lesssim \|u_1\|^{\theta_1\theta_2}_{L^2}\|\Lambda^s_1u_1\|^{\theta_1(1-\theta_2)}_{L^2}\|u_1\|^{1-\theta_1}_{H^k}\|\Lambda^{\sigma}u_2\|_{L^2},
\end{split}
\end{equation*}
and
\begin{equation*}
\begin{split}
\|\Lambda^{\frac12-\sigma_1}_1u_2^\sharp\|_{L^2_{x_1}L^{\infty}_{x_2}}\|\Lambda^{\sigma}u_1^\sharp\|_{L^2}&\lesssim \|\Lambda^{\frac12-\sigma_1}_1\Lambda^{\sigma}u_2\|^{\theta_3}_{L^2}\|\Lambda^{\frac12-\sigma_1}_1\Lambda^{k-\frac12+\sigma_1}_2u_2\|^{1-\theta_3}_{L^2}\|\Lambda^{\sigma}u_1\|_{L^2}\\&\lesssim  \|\Lambda^{\sigma}u_2\|^{\theta_3\theta_4}_{L^2}\|\Lambda^s_1\Lambda^{\sigma}u_2\|^{\theta_3(1-\theta_4)}_{L^2}\|u_2\|^{1-\theta_3}_{H^k}\|u\|^{\dot{\theta}_5}_{L^2}\|u\|^{1-\dot{\theta}_5}_{H^k}.
\end{split}
\end{equation*}
Here $\theta_1,\theta_2, \theta_3, \theta_4$  are given as in \eqref{11} and \eqref{2222}. 
While, 
$$
\dot{\theta}_5=\frac{k-\sigma}{k}.
$$
It's obvious that 
$
\dot{\theta}_5>\theta_5,
$
the same integrability argument as in
\eqref{eq2} and \eqref{eq3} yields that 
\begin{equation}\label{4.23}
\begin{split}
&\int_0^t|Q_2+Q_3|\,\mathrm{d}\tau\\&\lesssim E^{\frac{\theta_1}{2}}_1(t)E^{\frac{1-\theta_1}{2}}_2(t)E^{\frac12}_5(t)  E^{\frac12}_4(t)\left(\int_0^t(1+\tau)^{-\frac{\sigma_1\theta_1+\sigma_1+1}{s(2-\theta_1(1-\theta_2))}}\,\mathrm{d}\tau\right)^{\frac{2-\theta_1(1-\theta_2)}{2}}\\&\,\,\quad+E^{\frac{\theta_3}{2}}_5(t)E^{\frac{2-\theta_3-\dot{\theta}_5}{2}}_2(t)E^{\frac{\dot{\theta}_5}{2}}_1(t)E^{\frac12}_4(t)\left(\int_0^t(1+\tau)^{-\frac{(\sigma_1+1)\theta_3+\sigma_1\theta_5}{s(2-\theta_3(1-\theta_4))}}\,\mathrm{d}\tau\right)^{\frac{2-\theta_3(1-\theta_4)}{2}}
\\&\lesssim E^{\frac{\theta_1}{2}}_1(t)E^{\frac{1-\theta_1}{2}}_2(t)E^{\frac12}_5(t)  E^{\frac12}_4(t)+E^{\frac{\theta_3}{2}}_5(t)E^{\frac{2-\theta_3-\dot{\theta}_5}{2}}_2(t)E^{\frac{\dot{\theta}_5}{2}}_1(t)E^{\frac12}_4(t) \\&\lesssim E^{\frac32}_1(t)+E^{\frac32}_2(t)+E^{\frac32}_4(t)+E^{\frac32}_5(t).
\end{split}
\end{equation}
{\bf{Estimate of $Q_4$:}}
Using $\partial_2u_2=-\partial_1u_1$ and handling like for $Q_1$, we then derive 
\begin{equation}\label{4.4}
\begin{split}
\int_0^t|Q_4|\,\mathrm{d}\tau
&\lesssim  \int_0^t\|\Lambda^{-\sigma_1}_1\Lambda^{\sigma-1}(\partial_1u_1u_2)\|_{L^2}\|\Lambda_1^{-\sigma_1-1}\Lambda^{\sigma}u_2\|_{L^2}\,\mathrm{d}\tau 
 \\&\lesssim  \int_0^t \big(\|\partial_1u_1\|_{L^2}\|u_2\|_{L^2}+\|\partial_1u\|_{L^2}^{\frac32-\frac1{2s}}
\|\Lambda_1^{1+s}u\|_{L^2}^{\frac1{2s}}
\|u_2\|_{L^2}^{\frac12}\\&\,\,\quad+\|\partial_1u_1\|_{L^2}\|u_2\|^{\frac12}_{L^2}\|\partial_2u_2\|^{\frac12}_{L^2}\\&\,\,\quad+\|\partial_1u_1\|_{L^2}\|u_2\|^{\frac12}_{L^2}\|\partial_1u_2\|^{\frac12}_{L^2}\big)\|\Lambda_1^{-\sigma_1-1}\Lambda^{\sigma}u_2\|_{L^2}\,\mathrm{d}\tau  \\&\lesssim
E_4^{\frac12}(t)E_3^{\frac12}(t)E_1^{\frac12}(t)\int_0^t (1+\tau)^{-\frac{2\sigma_1+1}{2s}}\,\mathrm{d}\tau
\\&\,\,\quad+
E_4^{\frac12}(t)E_3^{\frac34}(t)E_1^{\frac14}(t)\left(\int_0^t(1+\tau)^{-\frac{4\sigma_1+3}{4s-1}}\,\mathrm{d}\tau\right)^{\frac{4s-1}{4s}}\\&\,\,\quad+
E_4^{\frac12}(t)E_3^{\frac34}(t)E_1^{\frac14}(t)\int_0^t(1+\tau)^{-\frac{4\sigma_1+3}{4s}} \,\mathrm{d}\tau
\\&\lesssim E^{\frac32}_1(t)+E^{\frac32}_3(t)+E^{\frac32}_4(t).
\end{split} 
\end{equation}

In the second inequality of \eqref{4.4}, the factor 
$\|\Lambda_1^{1+s}u\|_{L^2}^{\frac{1}{2s}}$ uses $\frac{1}{4s}$ of the time
H\"older's exponent. The remaining exponent is $\frac{4s-1}{4s}$, and
$\frac{4\sigma_1+3}{4s-1}>1$, so the displayed time integral is uniformly
bounded.

Combining the inequalities \eqref{4.1}, \eqref{4.23} and \eqref{4.4} together, we conclude the proof of Lemma \ref{lem4}.
\end{proof}

Lemma \ref{lem4} has controlled the negative-order norms of $\Lambda^{\sigma}u_2$, we now turn to deal with its time-weighted decay estimate.
\begin{lemma}\label{lem3.5}
Assume that 
$$
k>\max\left\{\frac{3\sigma-\sigma_1-4\sigma_1^2}
{2(\sigma-\sigma_1)},\frac{4\sigma^2_1-5\sigma_1+2}{2(1-2\sigma_1)},
\frac{1+\sigma_1-2\sigma+2\sigma\sigma_1}{3-6\sigma_1},1
\right\},
$$
and the energies are defined as in \eqref{eq1}, then we have
$$
E_5(t)\lesssim E_5(0)+E^{\frac{s}{1+\sigma_1}}_4(t)E^{\frac{1+\sigma_1-s}{1+\sigma_1}}_5(t)+E^{\frac32}_1(t)+E^{\frac32}_2(t)+E^{\frac32}_3(t)+E^{\frac32}_5(t).
$$
\end{lemma}
\begin{proof}
Applying $\Lambda^\sigma$ to \eqref{u_2}, taking the $L^2$ inner product
with $\Lambda^\sigma u_2$, and multiplying by the time weight
$(1+t)^{\frac{\sigma_1+1}{s}}$, we obtain
\begin{equation*}
\frac12\frac{\mathrm{d}}{\mathrm{d}t}\left((1+t)^{\frac{\sigma_1+1}{s}}\|\Lambda^{\sigma}u_2\|^2_{L^2}\right)+\nu(1+t)^{\frac{\sigma_1+1}{s}}\|\Lambda^{s}_1\Lambda^{\sigma}u_2\|^2_{L^2}= M,
\end{equation*}
in which 
\begin{equation*}
\begin{split}
M&= \frac12\frac{\sigma_1+1}{s}(1+t)^{\frac{\sigma_1+1-s}{s}}\|\Lambda^{\sigma}u_2\|^2_{L^2} +(1+t)^{\frac{\sigma_1+1}{s}}
\langle\Lambda^{\sigma}[(-\Delta)^{-1}\partial_1\nabla^\bot \cdot (u\cdot\nabla u)],\Lambda^{\sigma}u_2\rangle
\\&= \frac12\frac{\sigma_1+1}{s}(1+t)^{\frac{\sigma_1+1-s}{s}}\|\Lambda^{\sigma}u_2\|^2_{L^2} 
-(1+t)^{\frac{\sigma_1+1}{s}}<\Lambda^{\sigma}[(-\Delta)^{-1}\partial_1\partial_2\partial_1(u_1u_1)],\Lambda^{\sigma}u_2>\\&\,\,\quad-(1+t)^{\frac{\sigma_1+1}{s}}<\Lambda^{\sigma}[(-\Delta)^{-1}\partial_1\partial_2\partial_2(u_2u_1)],\Lambda^{\sigma}u_2>\\&\,\,\quad+(1+t)^{\frac{\sigma_1+1}{s}}<\Lambda^{\sigma}[(-\Delta)^{-1}\partial_1\partial_1\partial_1(u_1u_2)],\Lambda^{\sigma}u_2>\\&\,\,\quad+(1+t)^{\frac{\sigma_1+1}{s}}<\Lambda^{\sigma}[(-\Delta)^{-1}\partial_1\partial_1\partial_2(u_2u_2)],\Lambda^{\sigma}u_2> \\&\triangleq M_1+M_2+M_3+M_4+M_5.
\end{split}
\end{equation*}
{\bf{Estimate of $M_1$:}} Interpolation inequality  directly yields
\begin{equation}\label{3.51}
\begin{split}
\int_0^t|M_1|\,\mathrm{d}\tau& \lesssim\int_0^t(1+\tau)^{\frac{\sigma_1+1-s}{s}}\|\Lambda^{-\sigma_1-1+s}_1\Lambda^{\sigma}u_2\|^{\frac{2s}{1+\sigma_1}}_{L^2}\|\Lambda^{s}_1\Lambda^{\sigma} u_2\|^{\frac{2(1+\sigma_1-s)}{1+\sigma_1}}\,\mathrm{d}\tau\\&\lesssim \int_0^t\|\Lambda^{-\sigma_1-1+s}_1\Lambda^{\sigma}u_2\|^{\frac{2s}{1+\sigma_1}}_{L^2}\left((1+\tau)^{\frac{\sigma_1+1}{2s}}\|\Lambda^{s}_1\Lambda^{\sigma}u_2\|_{L^2}\right)^{\frac{2(1+\sigma_1-s)}{1+\sigma_1}}\,\mathrm{d}\tau\\&\lesssim E^{\frac{s}{1+\sigma_1}}_4(t)E^{\frac{1+\sigma_1-s}{1+\sigma_1}}_5(t).
\end{split}
\end{equation}
{\bf{Estimate of $M_2$:}}
Since {$0<s-\sigma<\sigma<\frac{1}{2}$}, we may choose
{$q_2=\frac{2}{1+2(s-\sigma)}\in(1,2)$} such that
$$
s-\sigma=\frac{2-q_2}{2q_2}.
$$
By the boundness of Calderón--Zygmund operators, we can derive 
\begin{equation*}
\begin{split}
|M_2|&=\big|(1+t)^{\frac{\sigma_1+1}{s}}<\Lambda^{\sigma}[(-\Delta)^{-1}\partial_1\partial_2\partial_1(u_1u_1)],\Lambda^{\sigma}u_2>\big|\\&\lesssim (1+t)^{\frac{\sigma_1+1}{s}}\|\Lambda^{\sigma-1}\Lambda^{1-s}_1(\partial_1u_1u_1)\|_{L^2}\|\Lambda^s_1\Lambda^{\sigma}u_2\|_{L^2}.
\end{split}
\end{equation*}
The first term on the right-hand side can be estimated like 
\begin{equation}\label{ri}
\begin{split}
&\|\Lambda^{\sigma-1}\Lambda^{1-s}_1(\partial_1u_1u_1)\|_{L^2}\\&\lesssim \|\Lambda^{s-1}_1\Lambda^{\sigma-s}_2\Lambda^{1-s}_1(\partial_1u_1u_1)\|_{L^2}\\&\lesssim \|\Lambda^{\sigma-s}_2(\partial_1u_1u_1)\|_{L^2}\\&\lesssim\|\partial_1u_1u_1\|_{L^2_{x_1}L^{q_2}_{x_2}}\\&\lesssim \|\partial_1u_1\|_{L^2}\|u_1\|_{L^{\infty}_{x_1}L^{\frac{2q_2}{2-q_2}}_{x_2}}\\&\lesssim \|\partial_1u_1\|_{L^2}\|(\Lambda^{\frac12-\widetilde\delta}_1+\Lambda^{\frac12+\widetilde\delta}_1)\Lambda^{\frac{2q_2-2}{2q_2}}_2u_1\|_{L^2} \\&\lesssim
 \|\partial_1u_1\|_{L^2}\|(\Lambda^{\frac12-\widetilde\delta}_1+\Lambda^{\frac12+\widetilde\delta}_1)u_1\|^{\Theta}_{L^2}\|(\Lambda^{\frac12-\widetilde\delta}_1+\Lambda^{\frac12+\widetilde\delta}_1)\Lambda^{k-\frac12-\widetilde\delta}_2u_1\|^{1-\Theta}_{L^2}\\&\lesssim\|\partial_1u_1\|_{L^2}\|(\Lambda^{\frac12-\widetilde\delta}_1+\Lambda^{\frac12+\widetilde\delta}_1)u_1\|^{\Theta}_{L^2}\|u\|^{1-\Theta}_{H^k},
\end{split}
\end{equation}
where
$$\Theta=\frac{k+s-\sigma-1-\widetilde\delta}
{k-\frac12-\widetilde\delta}.$$
In particular,

\begin{equation*}
\int_0^t\|\Lambda^{\frac12-\widetilde\delta}_1u_1\|^{2\Theta}_{L^2}\,\mathrm{d}\tau\lesssim \int_0^t \|u_1\|^{(\frac12+\widetilde\delta)2\Theta}_{L^2}\|\partial_1u_1\|^{(\frac12-\widetilde\delta)2\Theta}_{L^2} \,\mathrm{d}\tau,
\end{equation*}
\begin{equation*}
\int_0^t\|\Lambda^{\frac12+\widetilde\delta}_1u_1\|^{2\Theta}_{L^2}\,\mathrm{d}\tau\lesssim \int_0^t \|u_1\|^{(\frac12-\widetilde\delta)2\Theta}_{L^2}\|\partial_1u_1\|^{(\frac12+\widetilde\delta)2\Theta}_{L^2} \,\mathrm{d}\tau.
\end{equation*}
We choose $\widetilde\delta=\frac{1}{2}-\sigma>0$. The horizontal
$L^\infty$ estimate in \eqref{ri} turns to be low--high frequency estimate,
not the endpoint embedding of $\dot H^{\frac{1}{2}}$.
 Applying Cauchy--Schwarz's inequality on
$|\xi_1|\le1$ by $|\xi_1|^{-\sigma}$, and on $|\xi_1|>1$ by
$|\xi_1|^{-(1-\sigma)}$. Both squared weights are integrable in their
respective regions. In particular, 
\[
\Theta=\frac{k+s-\frac32}{k+\sigma-1},\qquad 0<\Theta<1.
\]
The two horizontal interpolation estimates above therefore give

\begin{equation}\label{3.211}
\int_0^t\|(\Lambda^{\frac12-\widetilde\delta}_1+\Lambda^{\frac12+\widetilde\delta}_1)u_1\|^{2\Theta}_{L^2}\,\mathrm{d}\tau\lesssim \big(E^{\Theta}_1(t)+E^{\Theta}_3(t)\big)\int_0^t (1+\tau)^{-\frac{\sigma_1+\sigma}{s}\Theta}\,\mathrm{d}\tau.
\end{equation}
Moreover, the assumptions on $k$ together with $s<2\sigma_1$ imply 
$$
k>\frac{\frac32(\sigma_1+\sigma)-s(1+\sigma_1)}{\sigma_1+\sigma-s},
$$
which ensures, 
$$
\frac{\sigma_1+\sigma}{s}\Theta>1.
$$
Therefore by \eqref{ri} and \eqref{3.211}  we derive that 
\begin{equation}\label{3.52}
\begin{split}
\int_0^t|M_2|\,\mathrm{d}\tau &\lesssim \sup_{0\le \tau\le t} (1+\tau)^{\frac{\sigma_1+1}{2s}}\|\partial_1u_1\|_{L^2}\left(\int_0^t (1+\tau)^{\frac{\sigma_1+1}{s}}\|\Lambda^s_1\Lambda^{\sigma}u_2\|^2_{L^2}\,\mathrm{d}\tau\right)^{\frac12}\\&\,\,\quad\cdot \sup_{0\le \tau\le t}\|u_1\|^{1-\Theta}_{H^k}\left(\int_0^t\|(\Lambda^{\frac12-\widetilde\delta}_1+\Lambda^{\frac12+\widetilde\delta}_1)u_1\|^{2\Theta}_{L^2}\,\mathrm{d}\tau\right)^{\frac12}\\&\lesssim  E^{\frac12}_{3}(t)E^{\frac12}_{5}(t)E^{\frac{1-\Theta}{2}}_2(t) \big(E^{\Theta}_1(t)+E^{\Theta}_3(t)\big)^{\frac{1}{2}}\\&\lesssim E^{\frac32}_1(t)+E^{\frac32}_2(t)+E^{\frac32}_3(t)+E^{\frac32}_5(t).
\end{split}
\end{equation}
{\bf{Estimate of $M_3+M_4$:}} {Calder\'on}--Zygmund {theory} and {Lemma \ref{lem:fractoinal-product-1}} yield that 
\begin{equation*}
\begin{split}
|M_3+M_4|&\lesssim (1+t)^{\frac{\sigma_1+1}{s}}\|\Lambda^{\sigma}\Lambda^{1-s}_1(u_1u_2)\|_{L^2}\|\Lambda^s_1\Lambda^{\sigma}u_2\|_{L^2}\\&\lesssim (1+t)^{\frac{\sigma_1+1}{s}}\Big(\|\Lambda^{\frac34-\frac s2}_1u_2^\sharp\|_{L^2_{x_1}L^{\infty}_{x_2}}\|\Lambda^{\frac34-\frac s2}_1\Lambda^{\sigma}u_1^\sharp\|_{L^2}\\&\,\quad+\|\Lambda^{\frac34-\frac s2}_1u_1^\sharp\|_{L^2_{x_1}L^{\infty}_{x_2}}\|\Lambda^{\frac34-\frac s2}_1\Lambda^{\sigma}u_2^\sharp\|_{L^2}\Big)\|\Lambda^s_1\Lambda^{\sigma}u_2\|_{L^2}.
\end{split}
\end{equation*}
We shall deal with the two product terms separately. 
For the first one, anisotropic type Sobolev
embedding and interpolation inequalities yield

\begin{equation*}
\begin{split}
 \|\Lambda^{\frac34-\frac s2}_1u_2^\sharp\|_{L^2_{x_1}L^{\infty}_{x_2}}\|\Lambda^{\frac34-\frac s2}_1\Lambda^{\sigma}u_1^\sharp\|_{L^2} &\lesssim \|\Lambda^{\frac34-\frac s2}_1\Lambda^{\sigma}_2u_2\|^{\Theta_1}_{L^2}\|\Lambda^{\frac34-\frac s2}_1\Lambda^k_2u_2\|^{1-\Theta_1}_{L^2}\|\Lambda^{\frac34-\frac s2}_1\Lambda^{\sigma}u_1\|_{L^2}\\&\lesssim \|\Lambda^{\sigma}u_2\|^{\Theta_1\Theta_2}_{L^2}\|\Lambda^s_1\Lambda^{\sigma}u_2\|^{\Theta_1(1-\Theta_2)}_{L^2}\|\Lambda^k_2u_2\|^{\Theta_3(1-\Theta_1)}_{L^2}\\&\,\,\quad\cdot
\|\Lambda^s_1\Lambda^k_2u_2\|^{(1-\Theta_3)(1-\Theta_1)}_{L^2}\|\Lambda^{\frac34-\frac s2}_1u_1\|^{\Theta_4}_{L^2}\|\Lambda^{\frac34-\frac s2}_1\Lambda^ku_1\|^{1-\Theta_4}_{L^2}\\&\lesssim  \|\Lambda^{\sigma}u_2\|^{\Theta_1\Theta_2}_{L^2}\|\Lambda^s_1\Lambda^{\sigma}u_2\|^{\Theta_1(1-\Theta_2)}_{L^2}\|\Lambda^k_2u_2\|^{\Theta_3(1-\Theta_1)}_{L^2}\\&\,\,\quad\cdot \|\Lambda^s_1\Lambda^k_2u_2\|^{(1-\Theta_3)(1-\Theta_1)}_{L^2}\|u_1\|^{\Theta_4\Theta_5}_{L^2}\|\partial_1u_1\|^{\Theta_4(1-\Theta_5)}_{L^2}\\&\,\,\quad\cdot \|\Lambda^ku_1\|^{\Theta_6(1-\Theta_4)}_{L^2}\|\Lambda^s_1\Lambda^k u_1\|^{(1-\Theta_6)(1-\Theta_4)}_{L^2}.
\end{split}
\end{equation*}
Here
$$
\Theta_1=\frac{k-\frac12}{k-\sigma},\qquad\Theta_2=\Theta_3=\Theta_6=\frac{\frac32s-\frac34}{s},
$$
$$
\Theta_4=\frac{k-\sigma}{k},\qquad\Theta_5=\frac14+\frac s2.
$$
Thus,
\begin{equation}\label{a}
\begin{split}
&\int_0^t(1+\tau)^{\frac{\sigma_1+1}{s}} \|\Lambda^{\frac34-\frac s2}_1u_2^\sharp\|_{L^2_{x_1}L^{\infty}_{x_2}}\|\Lambda^{\frac34-\frac s2}_1\Lambda^{\sigma}u_1^\sharp\|_{L^2}\|\Lambda^s_1\Lambda^{\sigma}u_2\|_{L^2}\,\mathrm{d}\tau\\&\lesssim \left(\sup_{0\le \tau\le t} (1+\tau)^{\frac{\sigma_1+1}{s}} \|\Lambda^{\sigma}u_2\|^2_{L^2}\right)^{\frac{\Theta_1\Theta_2}{2}}\left(\int_0^t(1+\tau)^{\frac{\sigma_1+1}{s}}\|\Lambda^s_1\Lambda^{\sigma}u_2\|^2_{L^2}\,\mathrm{d}\tau\right)^{\frac{\Theta_1(1-\Theta_2)+1}{2}}\\&\,\,\quad\cdot\sup_{0\le \tau\le t}\|u\|^{\Theta_3(1-\Theta_1)+\Theta_6(1-\Theta_4)}_{H^k}\left(\int_0^t \|\Lambda^s_1u\|^2_{H^k}\,\mathrm{d}\tau\right)^{\frac{(1-\Theta_3)(1-\Theta_1)+(1-\Theta_6)(1-\Theta_4)}{2}}\\&\,\,\quad\cdot\left(\sup_{0\le \tau\le t}(1+\tau)^{\frac{\sigma_1}{s}}\|u_1\|^2_{L^2}\right)^{\frac{\Theta_4\Theta_5}{2}}\left(\sup_{0\le \tau\le t}(1+\tau)^{\frac{\sigma_1+1}{s}}\|\partial_1u_1\|^2_{L^2}\right)^{\frac{\Theta_4(1-\Theta_5)}{2}}\\&\,\,\quad\cdot \left(\int_0^t\left(1+\tau\right)^{-\frac{(\sigma_1+1)\Theta_1+\sigma_1\Theta_4\Theta_5+(\sigma_1+1)\Theta_4(1-\Theta_5)-(\sigma_1+1)}{2s}z_1}\,\mathrm{d}\tau\right)^{\frac{1}{z_1}},
\end{split}
\end{equation}
where $z_1$ satisfies the following condition
$$
\frac{\Theta_1(1-\Theta_2)+1}{2}+\frac{(1-\Theta_3)(1-\Theta_1)+(1-\Theta_6)(1-\Theta_4)}{2}+\frac{1}{z_1}=1.
$$
The assumptions on $k$ together with $\frac{3}{4}<s<2\sigma_1$ imply 
$$
\begin{aligned}
&(-4s+2\sigma_1+3)k^2+(4s\sigma-2\sigma_1\sigma-\sigma_1-\sigma-1)k+2\sigma_1\sigma^2
\\
&>
(3-6\sigma_1)k^2-(1+\sigma_1-2\sigma+2\sigma\sigma_1)k>0 .
\end{aligned}
$$
Consequently,
$$
\frac{(\sigma_1+1)\Theta_1+\sigma_1\Theta_4\Theta_5+(\sigma_1+1)\Theta_4(1-\Theta_5)-(\sigma_1+1)}{2s}z_1>1,
$$
which ensures the desired time-weight integrability.

For the second product, using the similar interpolation argument we shall derive 
\begin{equation*}
\begin{split}
& \|\Lambda^{\frac34-\frac s2}_1u_1^\sharp\|_{L^2_{x_1}L^{\infty}_{x_2}}\|\Lambda^{\frac34-\frac s2}_1\Lambda^{\sigma}u_2^\sharp\|_{L^2}\\\lesssim &\|\Lambda^{\frac34-\frac s2}_1u_1\|^{\widetilde{\Theta}_1}_{L^2}\|\Lambda^{\frac34-\frac s2}_1\Lambda^{k-\frac34+\frac s2}_2u_1\|^{1-{\widetilde{\Theta}_1}}_{L^2}\|\Lambda^{\frac34-\frac s2}_1\Lambda^{\sigma}u_2\|_{L^2}\\\lesssim& \|u_1\|^{\widetilde{\Theta}_1\widetilde{\Theta}_2}_{L^2}\|\Lambda^s_1u_1\|^{\widetilde{\Theta}_1(1-\widetilde{\Theta}_2)}_{L^2}\|u\|^{1-\widetilde{\Theta}_1}_{H^k}\|\Lambda^{\sigma}u_2\|^{\widetilde{\Theta}_3}_{L^2}\|\Lambda^s_1\Lambda^{\sigma}u_2\|^{1-\widetilde{\Theta}_3}_{L^2}.
\end{split}
\end{equation*}
In which, 
$$
\widetilde{\Theta}_1=\frac{k+\frac s2-\frac54}{k+\frac s2-\frac34},\qquad \widetilde{\Theta}_2=\widetilde{\Theta}_3=\frac{\frac32s-\frac34}{s}.
$$
Consequently, we can obtain 
\begin{equation}\label{b}
\begin{split}
&\int_0^t(1+\tau)^{\frac{\sigma_1+1}{s}}\|\Lambda^{\frac34-\frac s2}_1u_1^\sharp\|_{L^2_{x_1}L^{\infty}_{x_2}}\|\Lambda^{\frac34-\frac s2}_1\Lambda^{\sigma}u_2^\sharp\|_{L^2}\|\Lambda^s_1\Lambda^{\sigma}u_2\|_{L^2}\,\mathrm{d}\tau\\&\lesssim \sup_{0\le \tau\le t} \left((1+\tau)^{\frac{\sigma_1}{s}}\|u_1\|^2_{L^2}\right)^{\frac{\widetilde{\Theta}_1\widetilde{\Theta}_2}{2}}\left(\int_0^t(1+\tau)^{\frac{\sigma_1}{s}}\|\Lambda^s_1u_1\|^2_{L^2}\,\mathrm{d}\tau\right)^{\frac{\widetilde{\Theta}_1(1-\widetilde{\Theta}_2)}{2}}\\&\,\,\cdot\sup_{0\le \tau\le t}\left((1+\tau)^{\frac{\sigma_1+1}{s}}\|\Lambda^{\sigma}u_2\|^2_{L^2}\right)^{\frac{\widetilde{\Theta}_3}{2}}\left(\int_0^t(1+\tau)^{\frac{\sigma_1+1}{s}}\|\Lambda^s_1\Lambda^{\sigma}u_2\|^2_{L^2}\,\mathrm{d}\tau\right)^{\frac{2-\widetilde{\Theta}_3}{2}}\\&\,\,\cdot\sup_{0\le \tau\le t} \|u\|^{1-\widetilde{\Theta}_1}_{H^k}\left(\int_0^t(1+\tau)^{-\frac{\sigma_1\widetilde{\Theta}_1}{2s}z_2}\,\mathrm{d}\tau\right)^{\frac{1}{z_2}},
\end{split}
\end{equation}
where $z_2$ satisfies the condition
$$
\frac{\widetilde{\Theta}_1(1-\widetilde{\Theta}_2)}{2}+\frac{2-\widetilde{\Theta}_3}{2}+\frac{1}{z_2}=1.
$$
The assumptions on $k$ together with $s<2\sigma_1$ then imply
$$
k>\frac{4s^2-2s\sigma_1-10s+5\sigma_1+6}{6+4\sigma_1-8s}.
$$
Which gives
$$
 \frac{\sigma_1\widetilde{\Theta}_1}{2s}z_2>1.
$$
Combing \eqref{a} with \eqref{b}, we obtain 
\begin{equation}\label{3.534}
\begin{split}
&\int_0^t|M_3+M_4|\,\mathrm{d}\tau \\&\lesssim E^{\frac{\Theta_1+1}{2}}_5(t)E^{\frac{2-\Theta_1-\Theta_4}{2}}_2(t)E^{\frac{\Theta_4\Theta_5}{2}}_1(t)E^{\frac{\Theta_4(1-\Theta_5)}{2}}_3(t)+E^{\frac{\widetilde{\Theta}_1}{2}}_1(t)E_5(t)E^{\frac{1-\widetilde{\Theta}_1}{2}}_2(t)\\&\lesssim E^{\frac32}_1(t)+E^{\frac32}_2(t)+E^{\frac32}_3(t)+E^{\frac32}_5(t).
\end{split}
\end{equation}
{\bf{Estimate of $M_5$:}} By {Calder\'on}--Zygmund {theory}, Lemmas \ref{2.3}, \ref{lem2.4} and Sobolev embedding theorem, we shall derive 
\begin{equation*}
\begin{split}
|M_5|&\lesssim (1+t)^{\frac{\sigma_1+1}{s}}\|\Lambda^{\sigma-1}\Lambda^{2-s}_1(u_2u_2)\|_{L^2}\|\Lambda^s_1\Lambda^{\sigma}u_2\|_{L^2}\\&\lesssim (1+t)^{\frac{\sigma_1+1}{s}}\|\Lambda^{2-s}_1(u_2u_2)\|_{L^{\frac{2}{2-\sigma}}}\|\Lambda^s_1\Lambda^{\sigma}u_2\|_{L^2}
\\&\lesssim  (1+t)^{\frac{\sigma_1+1}{s}} \|u_2\|_{L^{\frac{2}{1-\sigma}}}\|\Lambda^{2-s}_1u_2\|_{L^2}\|\Lambda^s_1\Lambda^{\sigma}u_2\|_{L^2}\\&\lesssim (1+t)^{\frac{\sigma_1+1}{s}}\|\Lambda^{\sigma}u_2\|_{L^2}\|\Lambda^{2-s}_1u_2\|_{L^2}\|\Lambda^s_1\Lambda^{\sigma}u_2\|_{L^2}\\&\lesssim (1+t)^{\frac{\sigma_1+1}{s}}\|\Lambda^{\sigma}u_2\|_{L^2}\|\Lambda^s_1u_2\|_{H^k}\|\Lambda^s_1\Lambda^{\sigma}u_2\|_{L^2}.
\end{split}
\end{equation*}
Which implies that 
\begin{equation}\label{3.55}
\begin{split}
\int_0^t|M_5|\,\mathrm{d}\tau&\lesssim \sup_{0\le \tau\le t} (1+\tau)^{\frac{\sigma_1+1}{2s}}\|\Lambda^{\sigma}u_2\|_{L^2}\left(\int_0^t\|\Lambda^s_1u_2\|^2_{H^k}\,\mathrm{d}\tau\right)^{\frac12}\\&\,\,\quad\cdot \left(\int_0^t(1+\tau)^{\frac{\sigma_1+1}{s}}\|\Lambda^s_1\Lambda^{\sigma}u_2\|^2_{L^2}\,\mathrm{d}\tau\right)^{\frac12}\\&\lesssim E^{\frac12}_2(t)E_5(t)\\&\lesssim E^{\frac32}_2(t)+E^{\frac32}_5(t).  
\end{split}
\end{equation}
Combining all the above estimates, namely \eqref{3.51}, \eqref{3.52}, \eqref{3.534} and \eqref{3.55}, we then complete the proof of Lemma \ref{lem3.5}.

\end{proof}

\subsection{Proof of the Theorem \ref{thm1.1}}
We now combine the a priori estimates obtained in the previous lemmas and
complete the proof of Theorem \ref{thm1.1}. Under the restrictions
$$
\frac34<s<2\sigma_1,
\qquad
\sigma_1<\sigma<\frac12,
$$
the lower bound of $k$ in Theorem \ref{thm1.1} ensures all the
regularity requirements in Lemmas \ref{lem1}--\ref{lem3.5}.
Therefore, all the preceding estimates are applicable.

Define the total energy functional as
$$
E_{\mathrm{total}}(t)
=
E_0(t)+E_2(t)+E_3(t)+E_4(t)+E_5(t).
$$
Notice that $E_1(t)$ is not included in $E_{\mathrm{total}}(t)$.
By \eqref{E1-interpolation}, there is 
$$
E_1(t)\lesssim E_{\mathrm{total}}(t),
\qquad
E^{\frac32}_1(t)
\lesssim E^{\frac32}_{\mathrm{total}}(t).
$$
We next deal with the homogeneous mixed terms in Lemmas
\ref{lem3} and \ref{lem3.5}. Since
$$
0<\frac{s}{1+\sigma_1}<1,
$$
Young's inequality with a small constant directly gives
$$
C
E_0^{\frac{s}{1+\sigma_1}}(t)
E_3^{\frac{1+\sigma_1-s}{1+\sigma_1}}(t)
\leq
\frac12E_3(t)+C E_0(t).
$$
Similarly,
$$
C
E_4^{\frac{s}{1+\sigma_1}}(t)
E_5^{\frac{1+\sigma_1-s}{1+\sigma_1}}(t)
\leq
\frac12E_5(t)+C E_4(t).
$$
Here and below, the value of $C$ may change  line by line.
It follows from Lemmas \ref{lem1}, \ref{lem2}, and \ref{lem4}
that
$$
E_0(t)+E_2(t)+E_4(t)
\lesssim
E_{\mathrm{total}}(0)
+
E^{\frac32}_{\mathrm{total}}(t).
$$
After absorbing the preceding $\frac12E_3(t)$ term into the left-hand
side of estimate in Lemma \ref{lem3}, we easily get 
$$
E_3(t)
\lesssim
E_3(0)+E_0(t)+E^{\frac32}_{\mathrm{total}}(t)
\lesssim
E_{\mathrm{total}}(0)
+
E^{\frac32}_{\mathrm{total}}(t).
$$
In the same way, Lemmas \ref{lem4} and \ref{lem3.5} will give
$$
E_5(t)
\lesssim
E_5(0)+E_4(t)+E^{\frac32}_{\mathrm{total}}(t)
\lesssim
E_{\mathrm{total}}(0)
+
E^{\frac32}_{\mathrm{total}}(t).
$$
Combining these estimates together, we arrive at
\begin{equation}\label{444}
E_{\mathrm{total}}(t)
\leq
C_0E_{\mathrm{total}}(0)
+
C_1E^{\frac32}_{\mathrm{total}}(t),
\end{equation}
where $C_0,C_1>0$ are independent of $t$.

The initial assumption \eqref{initial} implies
$$
E_{\mathrm{total}}(0)\leq C_{\mathrm{in}}\epsilon^2.
$$
Indeed, the initial values $E_3(0)$ and $E_5(0)$ are controlled by
$\|u_0\|_{H^k}^2$, whereas those of $E_0(0)$ and $E_4(0)$ are explicitly
included in \eqref{initial}.
Set
$$
K=2C_0C_{\mathrm{in}}.
$$
By local existence and continuity theorem, the following bound 
\begin{equation}\label{4444}
E_{\mathrm{total}}(t)\leq K\epsilon^2
\end{equation}
holds on a sufficiently short time interval. Let $T^*$ be the maximal
time for which \eqref{4444} holds. For every $t\in[0,T^*]$,
\eqref{444} gives
$$
\begin{aligned}
E_{\mathrm{total}}(t)
&\leq
C_0C_{\mathrm{in}}\epsilon^2
+
C_1E^{\frac12}_{\mathrm{total}}(t)
E_{\mathrm{total}}(t)
\\
&\leq
\frac{K}{2}\epsilon^2
+
C_1K^{\frac12}\epsilon E_{\mathrm{total}}(t).
\end{aligned}
$$
Choose $\epsilon>0$ sufficiently small so that
$$
C_1K^{\frac12}\epsilon\leq\frac14.
$$
Then
$$
E_{\mathrm{total}}(t)
\leq
\frac{2K}{3}\epsilon^2
<
K\epsilon^2,
\qquad t\in[0,T^*].
$$
This strictly improves the bootstrap assumption. Hence the standard continuity 
argument gives $T^*=\infty$,
and therefore
$$
E_{\mathrm{total}}(t)\lesssim\epsilon^2,
\qquad t\geq0.
$$
Finally, \eqref{E1-interpolation} also gives
$$
E_1(t)\lesssim\epsilon^2.
$$
This completes the proof of Theorem \ref{thm1.1}.

\vskip .3in
\centerline{\bf Acknowledgments}

\vskip .1in
The third author is partially supported by the National Natural Science Foundation of China (NNSFC) under grant 12271160.
~\\
~\\
\noindent{\bf{Data availability:}} No data was used for the research described in the article.
~\\
~\\
\noindent{\bf{Conflict of interest:}} The authors declare no conflicts of interest.


\vskip .4in
\begin{thebibliography}{89}
\bibitem{Ada}
R. A. Adams and J. J. F. Fournier,
\textit{Sobolev Spaces},
second edition, Pure and Applied Mathematics (Amsterdam), \textbf{140},
Elsevier/Academic Press, Amsterdam, 2003.
\bibitem{Dan}
H. Bahouri, J.-Y. Chemin and R. Danchin,
\textit{Fourier Analysis and Nonlinear Partial Differential Equations},
Grundlehren der mathematischen Wissenschaften, \textbf{343},
Springer, Heidelberg, 2011.
\bibitem{Bie}
Q. Bie, H. Fang, S. Wang and Y. Zhou,
Stability and sharp decay for the 3D incompressible anisotropic Navier--Stokes equations with fractional horizontal dissipation,
\textit{J. Differential Equations},
\textbf{463} (2026), 114167.
\bibitem{Cao}
C. Cao and J. Wu,
Global regularity for the 2D MHD equations with mixed partial dissipation and magnetic diffusion,
\textit{Adv. Math.},
\textbf{226} (2011), no. 2, 1803--1822.
\bibitem{C}
K. Choi and I.-J. Jeong,
Infinite growth in vorticity gradient of compactly supported planar vorticity near Lamb dipole,
\textit{Nonlinear Anal. Real World Appl.},
\textbf{65} (2022), Paper No. 103470, 20 pp.
\bibitem{zhang}
W. Deng and P. Zhang,
Large time behavior of solutions to 3-D MHD system with initial data near equilibrium,
\textit{Arch. Ration. Mech. Anal.},
\textbf{230} (2018), 1017--1102.
\bibitem{D1}
S. Denisov,
Infinite superlinear growth of the gradient for the two-dimensional Euler equation,
\textit{Discrete Contin. Dyn. Syst.},
\textbf{23} (2009), 755--764.
\bibitem{D2}
S. Denisov,
Double-exponential growth of the vorticity gradient for the two dimensional Euler equation,
\textit{Proc. Amer. Math. Soc.},
\textbf{143} (2015), 1199--1210.
\bibitem{D3}
S. Denisov,
The sharp corner formation in 2D Euler dynamics of patches:
infinite double exponential rate of merging,
\textit{Arch. Ration. Mech. Anal.},
\textbf{215} (2015), 675--705.
\bibitem{DWXZ}
B. Q. Dong, J. Wu, X. Xu and N. Zhu,
Stability and exponential decay for the 2D anisotropic Navier--Stokes equations with horizontal dissipation,
\textit{J. Math. Fluid Mech.},
\textbf{23} (2021), no. 4, Paper No. 100, 11 pp.
\bibitem{D}
T. D. Drivas, T. M. Elgindi and I.-J. Jeong,
Twisting in Hamiltonian flows and perfect fluids,
\textit{Invent. Math.},
\textbf{238} (2024), no. 1, 331--370.
\bibitem{Du}
Y. Du, W. Yang and Y. Zhou,
On the exponential stability of a stratified flow to the 2D ideal MHD equations with damping,
\textit{SIAM J. Math. Anal.},
\textbf{51} (2019), 5077--5102.
\bibitem{zz1}
W. Feng, W. Wang and J. Wu,
Nonlinear stability for the 2D incompressible MHD system with fractional dissipation in the horizontal direction,
\textit{J. Evol. Equ.},
\textbf{23} (2023), 32.

\bibitem{LLL}
L. Grafakos and S. Oh, 
The Kato-Ponce inequality, \textit{Comm. Partial Diﬀerential Equations.}
\textbf{39} (2014), 1128-1157.

\bibitem{YZ}
I.-J. Jeong, Y. Yao and T. Zhou,
Superlinear gradient growth for 2D Euler equation without boundary,
\textit{arXiv:2507.15739}.
\bibitem{J}
R. H. Ji, J. Wu and W. Yang,
Stability and optimal decay for the 3D Navier--Stokes equations with horizontal dissipation,
\textit{J. Differential Equations},
\textbf{290} (2021), 57--77.
\bibitem{Jiang}
F. Jiang and S. Jiang,
Asymptotic behaviors of global solutions to the two-dimensional non-resistive MHD equations with large initial perturbations,
\textit{Adv. Math.},
\textbf{393} (2021), 108084.
\bibitem{KS}
A. Kiselev and V. \v{S}ver\'ak,
Small scale creation for solutions of the incompressible two dimensional Euler equation,
\textit{Ann. of Math.},
\textbf{180} (2013), 1205--1220.
\bibitem{lai1}
S. H. Lai, J. Wu and J. Zhang,
Stabilizing phenomenon for 2D anisotropic magnetohydrodynamic system near a background magnetic field,
\textit{SIAM J. Math. Anal.},
\textbf{53} (2021), 6073--6093.
\bibitem{lai2}
S. H. Lai, J. Wu and J. Zhang,
Stabilizing effect of magnetic field on the 2D ideal magnetohydrodynamic flow with mixed partial damping,
\textit{Calc. Var. Partial Differential Equations},
\textbf{61} (2022), 126.
\bibitem{lai3}
S. H. Lai, J. Wu, J. Zhang and X. Zhao,
Global stability and sharp decay estimates for 3D MHD equations with only vertical dissipation near a background magnetic field,
\textit{Adv. Math.},
\textbf{486} (2026), 110747.
\bibitem{zz}
J. Li, H. Wang and D. Zheng,
Stability and sharp decay for 3D incompressible MHD system with fractional horizontal dissipation and magnetic diffusion,
\textit{Z. Angew. Math. Phys.},
\textbf{74} (2023), 44.
\bibitem{li}
C. Li, J. Wu and X. Xu,
Smoothing and stabilization effects of magnetic field on electrically conducting fluids,
\textit{J. Differential Equations},
\textbf{276} (2021), 368--403.
\bibitem{Lin1}
H. Lin, R. Ji, J. Wu and L. Yan,
Stability of perturbations near a background magnetic field of the 2D incompressible MHD equations with mixed partial dissipation,
\textit{J. Funct. Anal.},
\textbf{279} (2020), 108519.
\bibitem{Lin2}
H. Lin, J. Wu and Y. Zhu,
Stability and large-time behavior of 3D incompressible MHD equations with partial dissipation near a background magnetic field,
\textit{Arch. Ration. Mech. Anal.},
\textbf{249} (2025), Paper No. 26.

\bibitem{MB}
A. Majda and A. Bertozzi,
\textit{Vorticity and Incompressible Flow},
Cambridge University Press, 2002.
\bibitem{Ren}
X. Ren, J. Wu, Z. Xiang and Z. Zhang,
Global existence and decay of smooth solution for the 2-D MHD equations without magnetic diffusion,
\textit{J. Funct. Anal.},
\textbf{267} (2014), 503--541.
\bibitem{S1}
M. E. Schonbek,
$L^2$ decay for weak solutions of the Navier--Stokes equations,
\textit{Arch. Ration. Mech. Anal.},
\textbf{88} (1985), 209--222.
\bibitem{S2}
M. E. Schonbek,
Lower bounds of rates of decay for solutions to the Navier--Stokes equations,
\textit{J. Amer. Math. Soc.},
\textbf{4} (1991), 423--449.
\bibitem{S3}
M. E. Schonbek and T. P. Schonbek,
On the boundedness and decay of moments of solutions to the Navier--Stokes equations,
\textit{Adv. Differential Equations},
\textbf{5} (2000), 861--898.
\bibitem{WWZ}
Z. B. Wang, J. Wu and N. Zhu,
Stability and decay for the 2D anisotropic Navier--Stokes equations with fractional horizontal dissipation on $\mathbb{R}^{2}$,
\textit{arXiv:2601.16110}.
\bibitem{W}
M. Wiegner,
Decay results for weak solutions to the Navier--Stokes equations on $\mathbb{R}^{n}$,
\textit{J. London Math. Soc.},
\textbf{35} (1987), 303--313.
\bibitem{Zhu}
J. Wu and Y. Zhu,
Global solutions of the 3D incompressible MHD system with mixed partial dissipation and magnetic diffusion near an equilibrium,
\textit{Adv. Math.},
\textbf{377} (2021), 107466.
\bibitem{Zhu1}
Y. Zhou and Y. Zhu,
Global classical solutions of 2D MHD system with only magnetic diffusion on periodic domain,
\textit{J. Math. Phys.},
\textbf{59} (2018), 081505.
\bibitem{Z}
A. Zlato\v{s},
Exponential growth of the vorticity gradient for the Euler equation on the torus,
\textit{Adv. Math.},
\textbf{268} (2015), 396--403.

\end{thebibliography}
\end{document}